\documentclass{article}

\usepackage[a4paper, total={6.5in, 8.5in}]{geometry}
\usepackage{graphicx}
\usepackage{authblk}
\usepackage{soul}
\usepackage{amsmath}
\usepackage[
backend=biber,
giveninits=true,
style=numeric-comp,
sorting=none,
maxcitenames=1,
natbib=true,
sortlocale=en_US,
url=false, 
doi=true,
eprint=false
]{biblatex}
\usepackage{geometry}
\usepackage{atbegshi}
\usepackage{float,lscape}
\usepackage{multirow}
\usepackage{blindtext}
\usepackage{amssymb}
\usepackage{ragged2e}
\usepackage{blindtext}
\usepackage{longtable}
\usepackage{comment}
\usepackage{graphicx}
\usepackage{xcolor}
\usepackage{caption}
\usepackage{amsmath}
\usepackage{bm}
\usepackage{changepage}

\usepackage{amsfonts}
\usepackage{tabularx}
\usepackage{lipsum}
\usepackage[title]{appendix}
\usepackage{mathtools}
\usepackage{circuitikz}
\usetikzlibrary{patterns}
\usepackage{hyperref}
\usepackage{caption}
\usepackage{booktabs}
\usepackage{array}
\usepackage{siunitx}

\hypersetup{
	colorlinks   = true, 
	urlcolor     = blue, 
	linkcolor    = blue, 
	citecolor   = red 
}
\providecommand{\keywords}[1]
{
  \small	
  \textbf{\textit{Keywords---}} #1
}

\title{Variational Physics-informed Neural Operator (VINO) for Solving Partial Differential Equations}

\author[1]{Mohammad Sadegh Eshaghi\thanks{Corresponding author. Email: eshaghi.khanghah@iop.uni-hannover.de}}
\author[3]{Cosmin Anitescu}
\author[3,4]{Manish Thombre}
\author[3,5]{Yizheng Wang}
\author[1,2]{Xiaoying Zhuang}
\author[3]{Timon Rabczuk\thanks{Corresponding author. Email: timon.rabczuk@uni-weimar.de}}

\affil[1]{Chair of Computational Science and Simulation Technology, Institute of Photonics, Department of Mathematics and Physics, Leibniz University Hannover, 30167 Hannover, Germany}
\affil[2]{Department of Geotechnical Engineering, College of Civil Engineering, Tongji University, Siping Road 1239, 200092 Shanghai, China}
\affil[3]{Institute of Structural Mechanics, Bauhaus-Universität Weimar, Germany}
\affil[4]{Department of Mechanical Engineering, Indian Institute of Technology Bombay, Mumbai, India}
\affil[5]{Department of Engineering Mechanics, Tsinghua University, Beijing, China}

\date{}

\begin{document}

\maketitle

\begin{abstract}
Solving partial differential equations (PDEs) is a required step in the simulation of natural and engineering systems. The associated computational costs significantly increase when exploring various scenarios, such as changes in initial or boundary conditions or different input configurations. This study proposes the Variational Physics-Informed Neural Operator (VINO), a deep learning method designed for solving PDEs by minimizing the energy formulation of PDEs. This framework can be trained without any labeled data, resulting in improved performance and accuracy compared to existing deep learning methods and conventional PDE solvers. By discretizing the domain into elements, the variational format allows VINO to overcome the key challenge in physics-informed neural operators, namely the efficient evaluation of the governing equations for computing the loss. Comparative results demonstrate VINO's superior performance, especially as the mesh resolution increases.  As a result, this study suggests a better way to incorporate physical laws into neural operators, opening a new approach for modeling and simulating nonlinear and complex processes in science and engineering.

\end{abstract}


\section{Introduction} \label{sec:Introduction}

Recent advancements in machine learning have shown potential in addressing Partial Differential Equations (PDEs) \cite{raissi2019physics, yu2018deep, samaniego2020energy, li2020multipole, lu2019deeponet, lu2019deeponet, li2024physics, wang2021learning, brunton2020machine} with diverse applicability in wide applications such as mechanics \cite{sadegh2024deepnetbeam, kaewnuratchadasorn2024physics}  fluid mechanics \cite{raissi2020hidden, tartakovsky2020physics}, heat transfer \cite{lu2024surrogate}, bio-engineering \cite{kissas2020machine, sahli2020physics}, materials \cite{lu2020extraction, chen2020physics}, and finance \cite{elbrachter2022dnn, han2018solving}. Among these, a notable development is the operator learning approach for solving PDEs. Unlike traditional neural networks, which are limited to learning mappings between fixed-dimensional inputs and outputs, neural operators are designed to learn mappings between function spaces \cite{li2020fourier, li2020multipole, lu2019deeponet, lu2021learning}. Neural operators are recognized as universal approximators for any continuous operator \cite{lu2021learning}, equipping them to approximate a wide range of solutions for various parametric PDE families. Specifically, the solution operator maps input functions, such as initial and boundary conditions, to output solution functions \cite{liu2022learning, yang2022generic, wen2023real, pmlr-v202-bonev23a}. 

Although the Neural Operators generally can capture complex multi-scale systems, they struggle with generalization and accurately approximating the true operator when trained only on coarse-resolution data. Furthermore, it remains challenging for these models to generalize effectively to new scenarios and conditions beyond those present in the training data. In many cases, generating data may be costly, unavailable, or limited to low-resolution observations \cite{hersbach2020era5}. The physics-informed versions of neural operators \cite{li2024physics, wang2021learning} have been introduced to resolve the shortcomings of data-driven approaches but the accuracy of current physics-informed neural operators is not sufficient unless some paired input-output data is used. Therefore, the limitation that comes from the data-driven aspect of methods and relies on the availability of data, either from existing numerical solvers or direct physical observations, still exists. 

We believe that the inability of existing neural operator methods to learn PDE operators purely based on physics is due to the reliance on the strong form of governing equations, which involves higher order derivatives and is satisfied just at the collocation points in the domain. In addition, one of the major technical challenges in order to use the physical loss is to efficiently compute the derivatives of the solution \cite{li2024physics}. Moreover,  to incorporate the weak form of PDEs, the challenge of calculating the integral is also added. Therefore, we propose the \textbf{Variational Physics-Informed Neural Operator (VINO)}, which incorporates the energy form of physical laws into the networks. This form allows us to discretize the domain and use the shape functions to approximate the solution and completely resolve the differentiation and integration challenges by calculating them analytically. As a result, VINO can learn operators purely based on physical data and with consistent hyperparameters, shows better performance as the mesh size increases, and exhibits a converging behavior in finer mesh, unlike other methods such as PINO. In other words, we introduced a better way to incorporate the physical laws into neural operators. The code and data supporting this work will be publicly available at \href{https://github.com/eshaghi-ms/VINO}{https://github.com/eshaghi-ms/VINO}.

\section{Results}
\begin{figure}[htbp]
    \centering
    \includegraphics[width=0.6\textwidth]{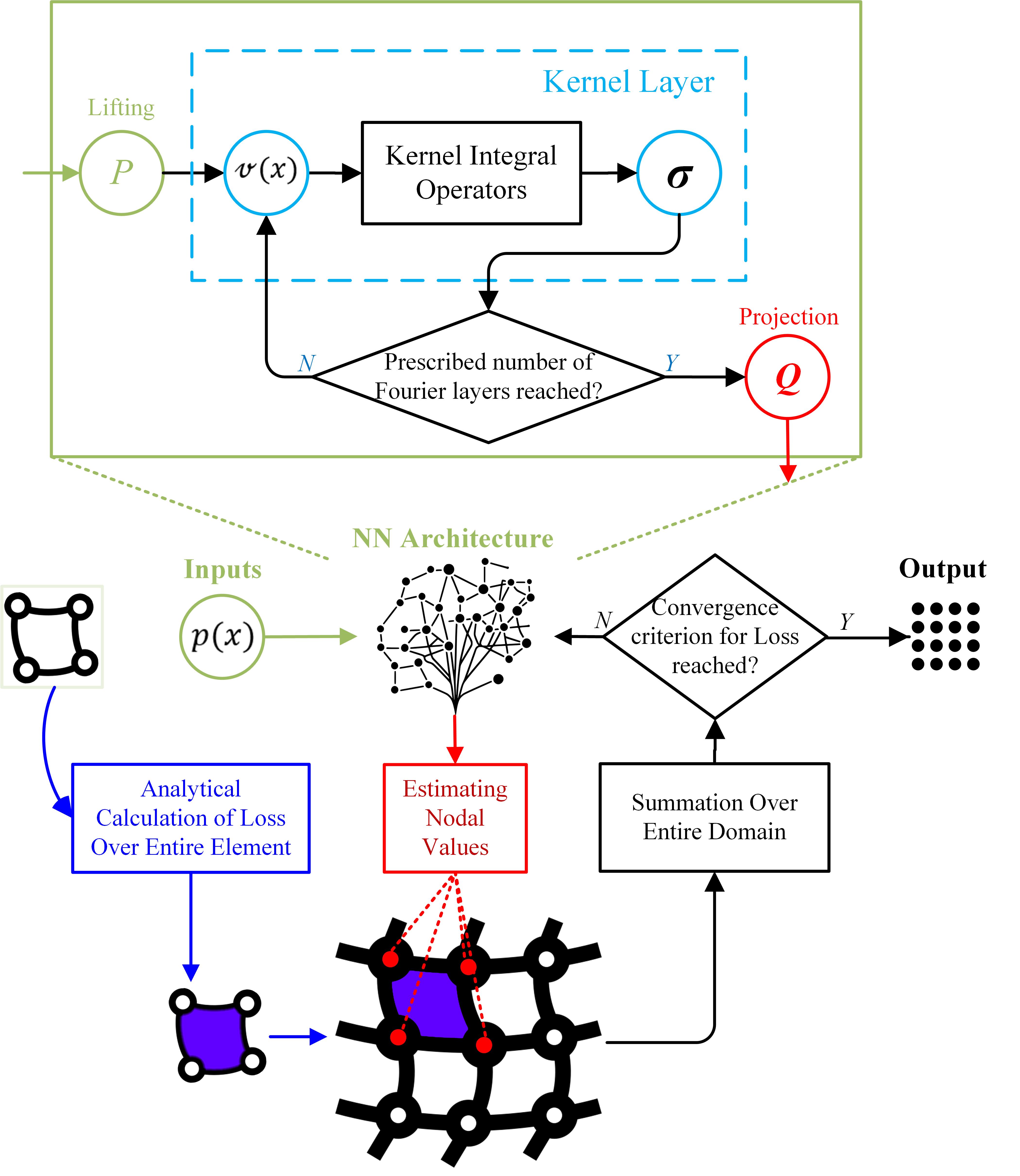}
    \caption{\textbf{Schematic Flowchart of VINO Method:} The flowchart outlines the VINO approach, which utilizes a variational loss function derived from the functional minimization of PDEs to learn the solution operator. This method eliminates the need for paired datasets, enabling greater robustness and adaptability across varying resolutions by employing domain discretization and deriving derivatives and integral analytically over each element.}
    \label{fig: VINO}
\end{figure}

The proposed Variational Physics-Informed Neural Operator framework is schematically depicted in Fig \ref{fig: VINO}. Inspired by the variational principle, VINO learns solution operators for PDEs without relying on paired datasets by minimizing the corresponding functional of PDEs as a loss function in a discretized domain. This discretization enables VINO to integrate and differentiate over each domain element analytically, removing dependency on limiting numerical methods or the computationally intensive automatic differentiation algorithm. 

\begin{figure}[htbp]
    \hspace{-50pt}
    \includegraphics[width=1.2\textwidth]{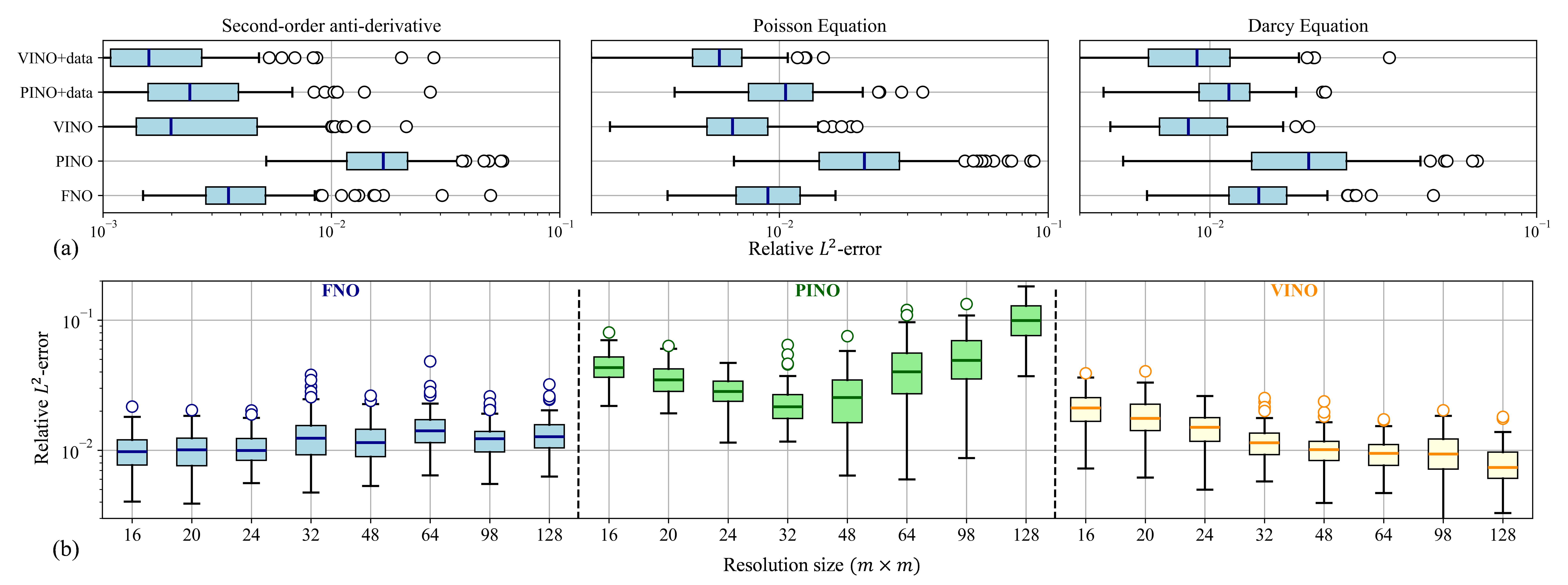}
    \caption{\textbf{Summary of the result.} (\textbf{a}) Comparison of relative \( L^2 \)-errors for three benchmark problems across different methods: FNO, PINO, VINO, their versions combined with data. Each box plot represents the distribution of errors for each method, showing the median, interquartile range, and outliers, indicating the robustness and accuracy of each method for the given problem. Lower \( L^2 \)-error values correspond to better performance, with VINO+data generally demonstrating lower error across all benchmarks. (\textbf{b}) Comparison of the relative \( L_2 \) error for FNO, PINO, and VINO methods on the Darcy equation test dataset across different mesh sizes. The boxplots show the distribution of errors for each method, the y-axis shows the \textbf{relative \( L_2 \) error} on a logarithmic scale, while the x-axis represents the \textbf{mesh size}, ranging from \( 16 \times 16 \) to \( 128 \times 128 \)}
    \label{fig: Summary}
\end{figure}

In the following, we illustrate the effectiveness of VINO over two other methods, Fourier Neural Operator (FNO) and Physics-Informed Neural Operator (PINO) on three benchmark examples and demonstrate its ability to solve more complicated cases across three further examples. A summary of methods' performance in the benchmark examples is presented in Table \ref{tab:Comparison}, and, graphically in Figure \ref{fig: Summary}a. In addition, we investigate the effect of mesh-grid size on the accuracy of different methods for solving the Darcy flow equation (for more detail refer to Section \ref{sec:darcy}) in Figure \ref{fig: Summary}b. All methods are trained with consistent hyperparameters across different mesh sizes and each boxplot illustrates the distribution of errors for the corresponding mesh size. \textbf{FNO} (blue boxplots) demonstrates relatively consistent performance with a slight upward trend across mesh sizes.  The median error shows a modest increase as the mesh size grows, indicating that FNO can effectively manage the challenges posed by finer meshes. In addition, the overall variability (as shown by the whiskers and outliers) remains relatively stable as the mesh resolution increases. For \textbf{PINO} (green boxplots), the relative \( L_2 \) error initially decreases with increasing mesh size, reaching its lowest median error around the \( 32 \times 32 \) mesh. However, beyond this mesh size, the error begins to rise, particularly at \( 98 \times 98 \) and \( 128 \times 128 \). Additionally, the variability increases with larger mesh sizes, as evidenced by a wider range and more noticeable outliers. The high error at smaller grid sizes can be due to small network parameters or insufficient accuracy in numerical integration and derivatives. As the grid size increases, PINO needs more iterations and more network parameters to reach the same accuracy, resulting in the \( 32 \times 32 \) mesh yielding the best performance with constant hyperparameters. For \textbf{VINO} (orange boxplots), the median relative error remains low overall, demonstrating a clear improvement in accuracy as the mesh size increases. The median relative \( L_2 \) error consistently decreases with finer mesh resolutions. Furthermore, variability in error diminishes for smaller meshes, indicating that VINO achieves more reliable accuracy improvements with increasing mesh size. Consequently, these results exhibit a converging behavior for VINO in finer mesh, with constant hyperparameters, reminiscent of the finite element method. 

\begin{table}[htbp]
    \centering
    \begin{tabular}{@{}p{3cm}p{2cm}p{2cm}p{2cm}p{2.2cm}p{2.2cm}@{}}
        \toprule
        \textbf{PDE Equation} & \textbf{FNO} & \textbf{PINO} & \textbf{VINO} & \textbf{PINO+data} & \textbf{VINO+data} \\
        \toprule
        Second-order & $0.54 \pm 0.61\%$ & $1.83 \pm 0.97\%$ & $0.36 \pm 0.36\%$ & $0.33 \pm 0.33\%$ & $0.25 \pm 0.35\%$ \\
        anti derivative & 94.19 s & 99.59 s & 101.04 s & 109.20 s & 112.89 s \\
        \midrule
        Poisson & $0.94 \pm 0.30\%$ & $2.63 \pm 1.92\%$ & $0.75 \pm 0.33\%$ & $1.11 \pm 0.52\%$ & $0.63 \pm 0.26\%$ \\
        & 220.9 s & 224.1 s & 236.5 s & 226.3 s & 236.7 s \\
        \midrule
        Darcy Flow & $1.49 \pm 0.58\%$ & $2.19 \pm 1.17\%$ & $0.93 \pm 0.31\%$ & $1.14 \pm 0.33\%$ & $0.98 \pm 0.46\%$ \\
        & 228.2 s & 251.2 s & 244.2 s & 251.3 s & 245.6 s \\
        \bottomrule
    \end{tabular}
    \caption{Statistical values of relative \(L^2\)-error and training time for different methods}
    \label{tab:Comparison}
\end{table}

\subsection{Comparative Numerical Results}

To illustrate the effectiveness of VINO, we present a series of numerical studies aimed at comparing with the best existing methods, including FNO and PINO. Because these methods are officially released using PyTorch frameworks, for fair comparison, we also developed our method in PyTorch for these comparative numerical experiments. However, we implemented the more efficient framework, JAX \cite{jax2018github} for solving some practical examples in the next section. 

In most cases, we represent random input functions \( p(x) \) using Gaussian random fields (GRF) \cite{williams2006gaussian} defined as 
\begin{equation}
    \begin{split}        
        p(x) \sim \mathcal{GP} \,(a; \, k_l(x_1, x_2)),
        \label{eq: GRF}
    \end{split}
\end{equation}
where \( k_l(x_1, x_2) \) is an exponential quadratic covariance kernel given by 
\begin{equation}
    \begin{split}        
        k_l(x_1, x_2) = \exp\left(-\frac{\|x_1 - x_2\|^2}{2l^2}\right),
        \label{eq: GRF_Kernel}
    \end{split}
\end{equation} 
and \(a\) and \( l \) represent the average of the field and positive length scale respectively. The parameter \( l \) is used to regulate the complexity of the sampled input functions, with larger values of \( l \) typically resulting in smoother \( p \). The numerical experiments were conducted on an NVIDIA A100-PCIE-40GB GPU and the networks are trained with mini-batch stochastic gradient descent, optimized using the Adam algorithm \cite{kingma2014adam}. 

\subsubsection{Second-order anti-derivative}

To demonstrate the capability of VINO in solving PDEs, let's first consider a simple ODE involving the second-order anti-derivative operator with the form
\begin{equation}
    \begin{split}        
        & \frac{d^2s(x)}{dx^2} + p(x) = 0, \quad x \in [0, 1], \\
        & s(0) = s(1) = 0.
        \label{eq: second-order anti derivative}
    \end{split}
\end{equation}

Given that \(p(x)\) is a known ODE parameter, the goal is to learn the operator \(\mathcal{G}\) mapping from \(p(x)\) to the solution \(s(x)\). To achieve this, we compare the results and performance of FNO, PINO, and VINO, with the operator denoted by \(\mathcal{G}_\theta\) with the same set of hyperparameters. Given a collection of paired parameters and solutions \(\{\mathbf{p}^{(i)}, \mathbf{s}^{(i)}\}\), where \(\mathbf{s}^{(i)} = \mathcal{G}({\mathbf{p}^{(i)}})\), the loss functions for FNO, \( \mathcal{L}_f\) is given by:
\begin{align}
    \begin{split}        
        \mathcal{L}_{f} & := \left\| \mathcal{G}-\mathcal{G}_{\theta} \right\| ^2_{L^2}
                        = \frac{1}{n} \sum_{i=1}^{n} \, \left\| \mathbf{s}^{(i)}
                        -\mathcal{G}_{\theta}(\mathbf{p}^{(i)}) \right\|^2_{L^2}, 
                        \label{eq: Problem-1_loss_data}
    \end{split}
\end{align}
and using the ODE itself, we have the loss function for PINO, \( \mathcal{L}_p\) as following:
\begin{align}
    \begin{split}        
        \mathcal{L}_{p} & := 
                        \frac{1}{n} \sum_{i=1}^{n} 
                        \left\| 
                        \frac{d^2}{dx^2} G_\theta(\mathbf{p}^{(i)})(x) + \mathbf{p}^{(i)}(x) 
                        \right\|^2_{L^2}
                        \label{eq: Problem-1_loss_pino}
    \end{split}
\end{align}
where \(\mathbf{p}^{(i)} = [p^{(i)}(x_1), p^{(i)}(x_2), \dots, p^{(i)}(x_m)]\) and the solutions \(\mathbf{s}^{(i)} = [s^{(i)}(x_1), s^{(i)}(x_2), \dots, s^{(i)}(x_m)]\), may be obtained by solving the ODE with a numerical method.

In addition, using the variational principle, a functional \(\Pi\) for the current problem can be defined in the form of
\begin{equation}
    \begin{split}        
        \Pi = & \int_\Omega \left[ \frac{1}{2} \left( \frac{ds}{dx} \right)^2 - p(x) \, s(x) \right] d\Omega
    \end{split}
\end{equation}
and the solution is a function which makes \(\Pi\) minimum so it also may be called the energy minimization principle. Therefore, after discretization and using approximative function \(\mathcal{G}_\theta^e(p) = \mathbf{N} \mathbf{\tilde{s}}_e\), which is discussed in details in Section \ref{sec: Method},  we can define the loss function for VINO, \( \mathcal{L}_v\) as:
\begin{equation}
    \begin{split}        
        \mathcal{L}_{v} := & \frac{1}{n} \sum_{i=1}^{n} \sum_e \left[ \frac{1}{2} \mathbf{\tilde{s}}_e^{(i) \, T} \int_{\Omega_e} \left(  \frac{d\mathbf{N}}{dx}^T\frac{d\mathbf{N}}{dx} \right) d\Omega_e \,\, \mathbf{\tilde{s}}_e^{(i)} - \mathbf{\tilde{s}}_e^{(i) \, T} \int_{\Omega_e} \left( \mathbf{N}^T \mathbf{N} \right) d\Omega_e \,\, \mathbf{\tilde{p}}_e^{(i)} \right] 
    \end{split}
\end{equation}
where \(\mathbf{N}\), \(\mathbf{\tilde{s}}_e\) and \(\mathbf{\tilde{p}}_e\) are the shape functions, nodal values of solution and parameter functions and \(e\) is the element index. The loss function also, in a simpler way, can be represented as follows: 
\begin{align}
    \begin{split}        
        \mathcal{L}_{v} & = \frac{1}{n} \sum_{i=1}^{n} \sum_e \left[ \frac{1}{2} \mathbf{\tilde{s}}_e^{(i) \, T} \mathbf{K}^e \, \mathbf{\tilde{s}_e}^{(i)} - \mathbf{\tilde{s}}_e^{(i) \, T} \mathbf{M}^e \, \mathbf{\tilde{p}}_e^{(i)} \right] 
        \label{eq: Problem-1_loss_vino}
    \end{split}
\end{align}
with element matrices \(\mathbf{K}^e\) and \(\mathbf{M}^e\) defined as:
\begin{equation}
    \begin{split}        
    \mathbf{K}^e &:= \int_{\Omega_e} \left(  \frac{d\mathbf{N}}{dx}^T\frac{d\mathbf{N}}{dx} \right) d\Omega_e
    \end{split}
\end{equation}
\begin{equation}
    \begin{split}        
    \mathbf{M}^e &:= \int_{\Omega_e} \left( \mathbf{N}^T \mathbf{N} \right) d\Omega_e
    \end{split}
\end{equation}
Therefore, on a uniform rectangular grid, the matrices \(K_e\) and \(M_e\) are the same for all elements which reduces the computational and memory cost.

We sample \(n = 1,000\) input functions \(p(x)\) from a Gaussian random field (GRF) with a length scale \(l = 0.1\). Additionally, we set \(m = 256\) sensors \(\{x_j\}_{j=1}^m\) uniformly spaced grid points in \([0, 1]\) and repeat the same procedure to generate a test dataset, which contains 100 different samples of random input functions and solutions. Finally, we evaluate and compare the performance of the three models using input functions that are not derived from the GRF, and therefore differ from both the training and test data. The corresponding solutions are obtained by analytically solving the ODE.

\begin{figure}[htbp]
    \centering
    \includegraphics[width=1.0\textwidth]{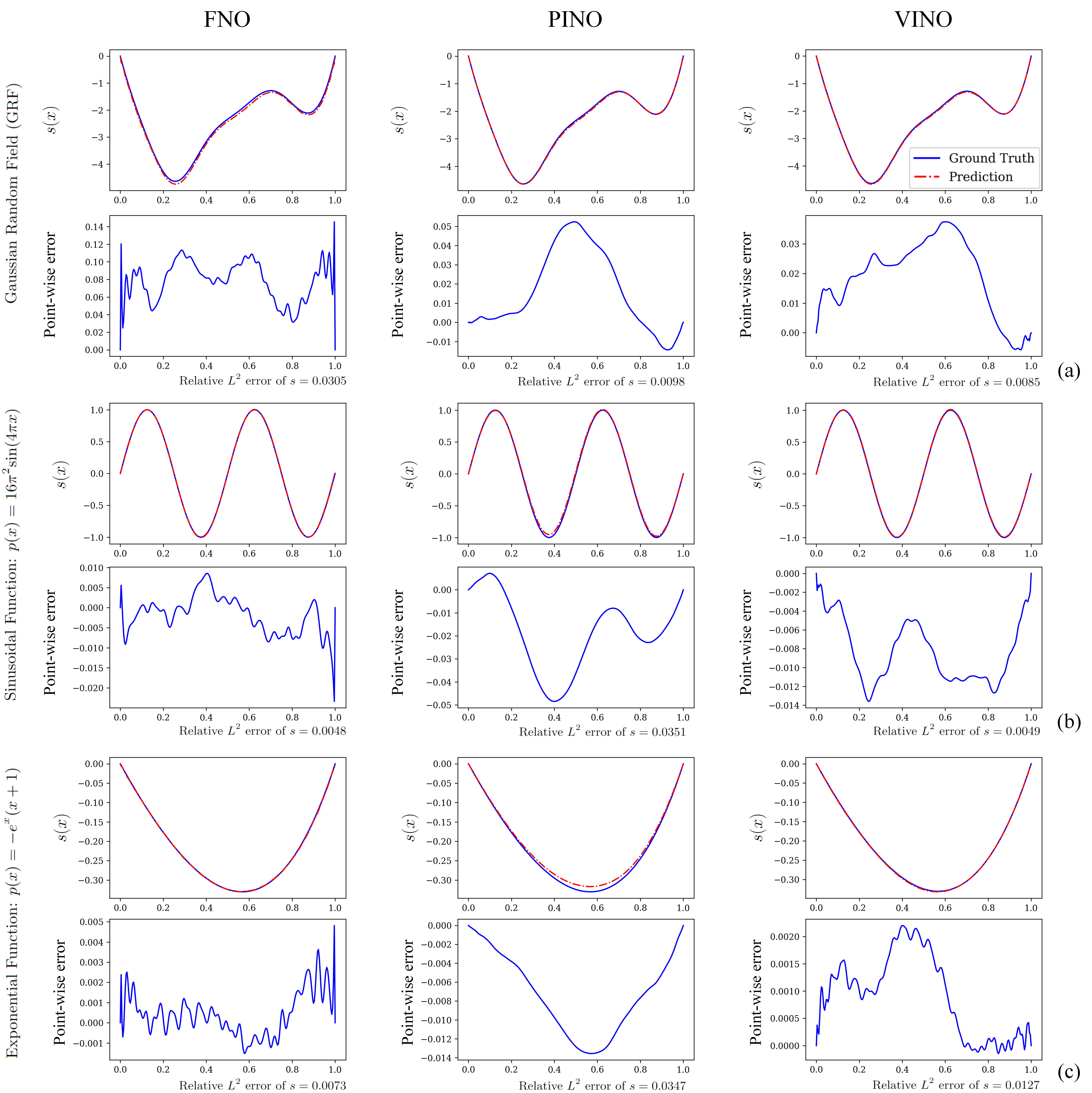}
    \vspace{-20pt}
    \caption{Comparison of the performance of FNO, PINO, and VINO for three representative test cases: (a) Gaussian Random Field (GRF), (b) Sinusoidal function, and (c) Exponential function. The top row in each case shows the predicted solutions (red) against the true solutions (blue), while the bottom row shows the point-wise error. The relative \(L^2\) errors are indicated for each model in each case.}
    \label{fig: Exp_01}
\end{figure}

We train the models by optimizing the specified loss function for 1,000 iterations using gradient descent with the Adam optimizer. Figure \ref{fig: Exp_01} presents the results of applying the FNO, PINO, and VINO models to different types of test functions. Specifically, the figure shows the predictions of each model (in red) against the true solutions (in blue) and their point-wise errors for three distinct test cases:

- Fig.\ref{fig: Exp_01}a - Gaussian Random Field (GRF): This row shows the performance of each model on a function sampled from the Gaussian random field, similar to the training data. The point-wise error below each plot reveals that while all models capture the overall trend of the solution, VINO exhibits a smaller error magnitude, as seen in the lower point-wise error values and the relative \(L^2\) errors.

- Fig.\ref{fig: Exp_01}b - Sinusoidal Function: In this case, the function is sinusoidal, differing from the GRF data used during training. The comparison between PINO and VINO reveals that VINO accurately captures both the phase and amplitude of the sine wave, while PINO struggles in some areas. Moreover, VINO shows a noticeably lower error compared to PINO, reflecting its improved ability to generalize.

- Fig.\ref{fig: Exp_01}c - Exponential Function: The exponential function further challenges the models. Here, the VINO model again shows a marked improvement over the PINO. Its point-wise error remains more contained, indicating that VINO generalizes better even for inputs that significantly differ from the training data.

In addition, we compare the performance of the PINO and VINO when combined with data on this parametric ODE problem. The first row in Table \ref{tab:Comparison} summarizes a statistical comparison of the relative $L_2$ prediction error for the output functions for different models, and their variants combined with data over 100 examples in the test dataset. 

We observe that VINO and its combination with data exhibit lower error means compared to PINO. In addition, VINO shows smaller deviations and variances compared to PINO, suggesting that they are more stable across different datasets. Therefore, the table and box plot indicate that VINO outperforms PINO. Note that the training time is approximately the same and VINO does not impose more computational demand. 

\subsubsection{Poisson Equation}

The next example involves two-dimensional Poisson equation, which is a second-order PDE with a source term \(p(x)\). The equation is given by:
\begin{equation}
    \begin{split}
        \nabla^2 s(x) + p(x) & = 0, \quad x \in \Omega = [0, 1]^2, \\  
        s(x) & = 0, \quad \text{in } \partial \Omega.
        \label{eq: poisson-2d}
    \end{split}
\end{equation}

Here,  the goal is to learn the operator \(\mathcal{G}\) mapping from PDE parameter \(p(x)\) to the solution \(s(x)\). Similar to the 1D case, we compare the performance of FNO, PINO, and VINO, denoted as \(\mathcal{G}_\theta\), using identical hyperparameter settings. For a collection of input-output pairs \(\{\mathbf{s}^{(i)}, \mathbf{p}^{(i)}\}\), where \(\mathbf{s}^{(i)} = \mathcal{G}(\mathbf{p}^{(i)})\), the loss functions for FNO is same as Eq. \ref{eq: Problem-1_loss_data}, but PINO and VINO loss function, with the same procedure, is defined as
\begin{align}
    \begin{split}
        \mathcal{L}_{p} & := \frac{1}{n} \sum_{i=1}^{n}
                        \left\| 
                        \frac{\partial^2}{\partial x^2} G_\theta(\mathbf{p}^{(i)}) + 
                        \frac{\partial^2}{\partial y^2} G_\theta(\mathbf{p}^{(i)}) + \mathbf{p}^{(i)} \right\|^2_{L^2}
                        \label{eq: Poisson2D_loss_pino}
    \end{split}
\end{align}
\begin{equation}
    \begin{split}        
        \mathcal{L}_{v} := & \frac{1}{n} \sum_{i=1}^{n} \sum_e \left[ \frac{1}{2} \mathbf{\tilde{s}}_e^{(i) \, T} \int_{\Omega_e} \left(  \nabla \mathbf{N}^T \, \nabla \mathbf{N} \right) d\Omega_e \,\, \mathbf{\tilde{s}}_e^{(i)} - \mathbf{\tilde{s}}_e^{(i) \, T} \int_{\Omega_e} \left( \mathbf{N}^T \mathbf{N} \right) d\Omega_e \,\, \mathbf{\tilde{p}}_e^{(i)} \right] 
        \label{eq: Poisson2D_loss_vino}
    \end{split}
\end{equation}

In this case, \(\mathbf{p}^{(i)} = [p^{(i)}(x_1), p^{(i)}(x_2), \dots, p^{(i)}(x_m)]\) and \(\mathbf{s}^{(i)} = [s^{(i)}(x_1), s^{(i)}(x_2), \dots, s^{(i)}(x_m)]\), with \(m = 64 \times 64\) grid points in the \([0, 1] \times [0, 1]\) domain. The input functions \(p(x)\) are sampled from a Gaussian random field (GRF). We use \(n = 1,000\) training samples and generate an additional test dataset of 100 different random input functions and solutions and train the models using gradient descent with the Adam optimizer

\begin{figure}[htbp]
    \centering
    \includegraphics[width=1.0\textwidth]{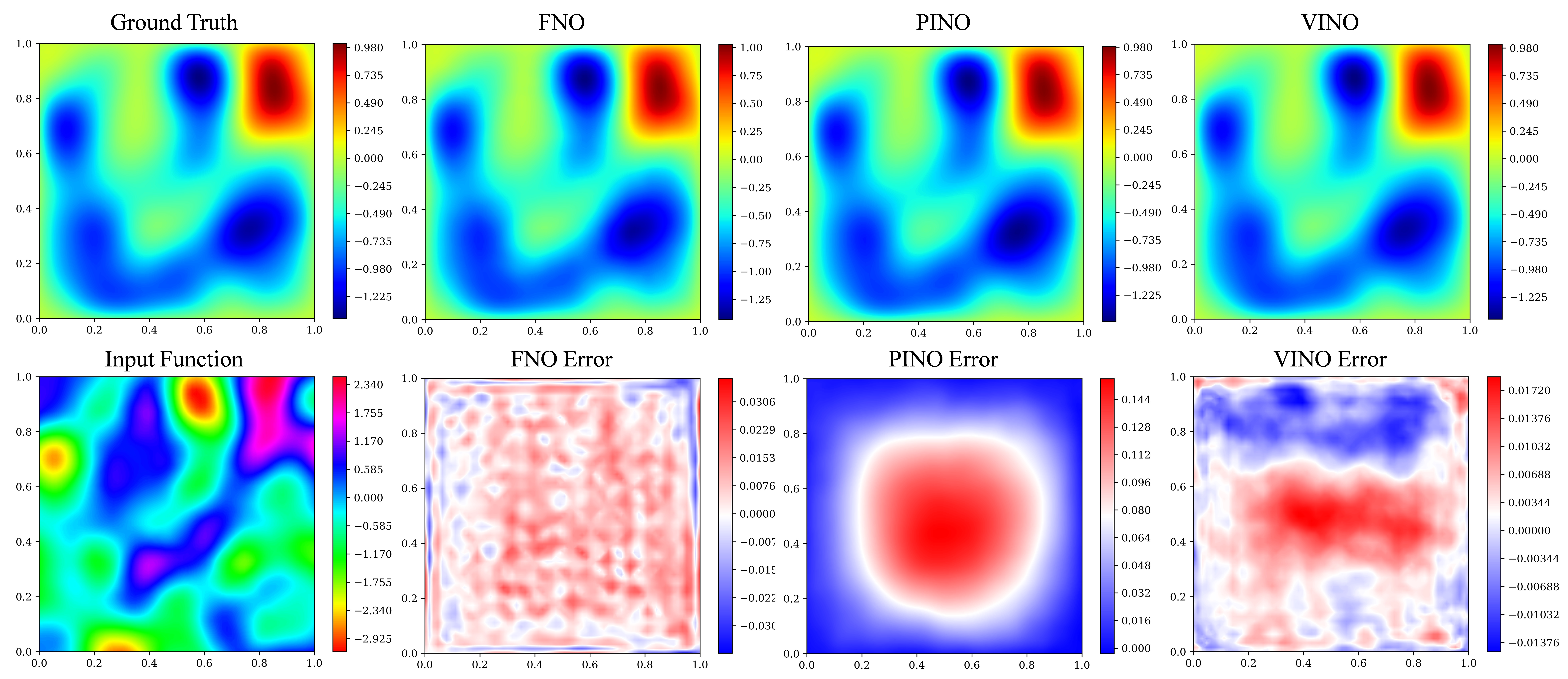}
    \vspace{-20pt}
    \caption{Performance comparison of FNO, PINO, and VINO for a sample from test dataset}
    \label{fig: Exp_02}
\end{figure}

Figure \ref{fig: Exp_02} illustrates the results for a sample from the test dataset. The figure shows that VINO performs better than PINO. In addition, as shown in Table \ref{tab:Comparison}, the second row, the mean errors highlight that VINO, both in its standalone and data-enhanced form, consistently achieves the lowest mean $L^2$ error of 0.00751 and 0.00631, respectively. This indicates that VINO is more effective at solving the Poisson equation than PINO.  In contrast, PINO not only has a higher mean error but also a greater deviation of 0.0192, indicating more variability in its performance across trials. Incorporating data improves the performance of both methods, with data-enhanced VINO achieving the lowest mean error and a reduced standard deviation of 0.0026. Notably, even with data augmentation, PINO's performance does not surpass that of VINO in its standalone form. As shown in the table, the training times and computational demands are very similar for all three methods.

\subsubsection{Darcy Flow} \label{sec:darcy}

As the next benchmark example, we focus on modeling two-dimensional subsurface flows through a porous medium with heterogeneous permeability fields. The high-fidelity synthetic simulation data are based on Darcy’s flow, explored in several neural operator studies \cite{li2020fourier, lu2022comprehensive, you2022nonlocal, you2022learning, li2020multipole}. The governing differential equation for this problem is given by:
\begin{equation}
    \begin{split}
        -\nabla \cdot (p(x) \nabla s(x)) = f(x), & \quad x \in [0,1]^2 \\  
        s(x) = 0, & \quad x \in \partial \, [0,1]^2
        \label{eq: Darcy-2d}
    \end{split}
\end{equation}

Here, \(p(x)\) is a piecewise constant conductivity field in \(L^\infty((0,1)^2; \mathbb{R}_+)\), and \(f(x) = 1\) is a constant source term function.   The goal is to learn the operator \(\mathcal{G}\) that maps the conductivity field \(p(x)\) to the solution, the hydraulic head, \(s(x)\), i.e., \(\mathcal{G}: p \mapsto s\). Despite the linearity of the PDE, the operator \(\mathcal{G}\) is nonlinear.

\begin{figure}[htbp]
    \centering
    \includegraphics[width=1.0\textwidth]{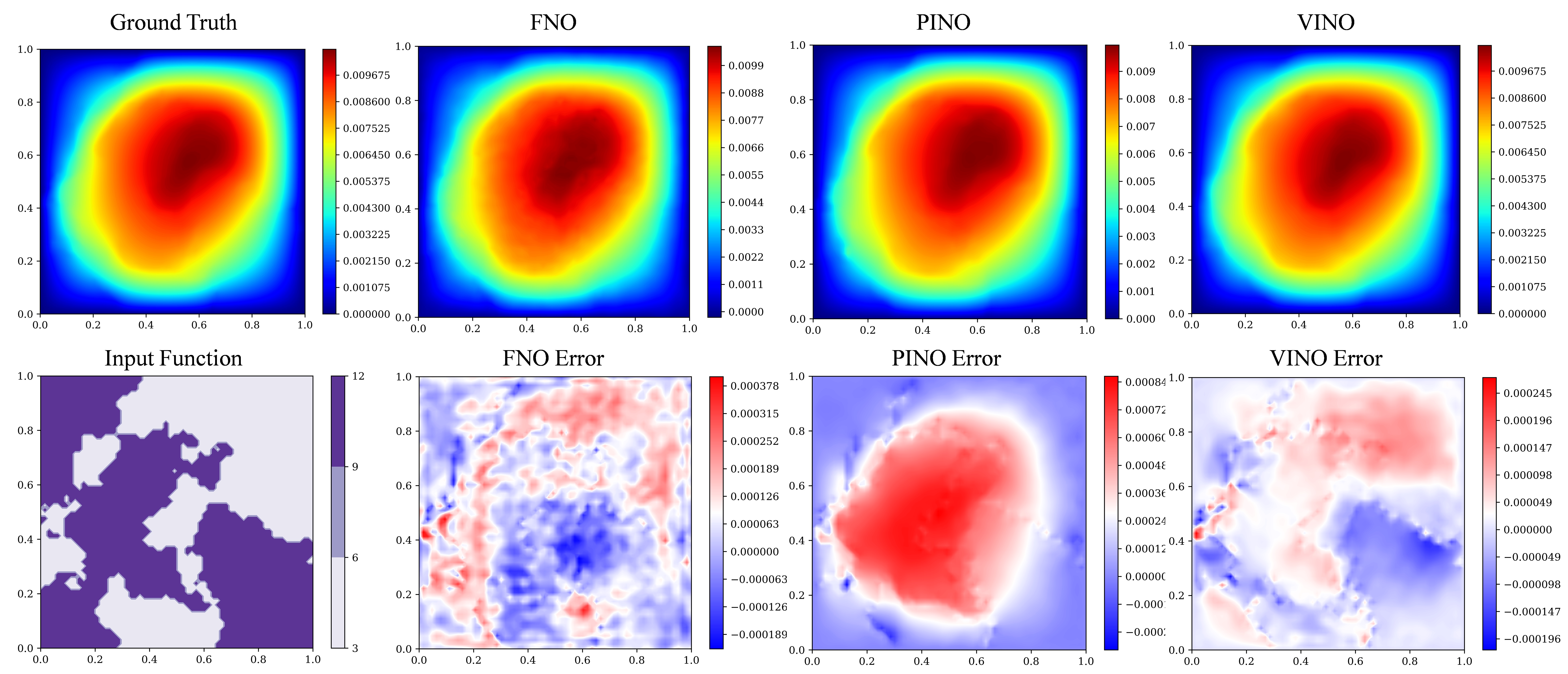}
    \caption{Performance comparison of FNO, PINO, and VINO for a sample from the test dataset in Darcy Equation.}
    \label{fig:Exp_03}
\end{figure}

As typical in subsurface flow simulations, the permeability \(p(x)\) is modeled as a piecewise constant function with random geometry, where the two values of \(p(x)\) have a ratio of 4. Specifically, 1100 samples of \(p(x)\) were generated according to the GRF distribution. Of these, 1000 samples were used for training, while the remaining were reserved for testing. 

As stated in Table \ref{tab:Comparison}, VINO has the lowest mean error, indicating the best overall performance, while PINO has the highest mean error and the largest standard deviation, reflecting more variability in its results. The addition of data reduces both the mean error and standard deviation for PINO and VINO, improving their performance and consistency. The training time across methods is relatively consistent, with small variations.

The box plot in Figure \ref{fig: Summary}a further illustrates the distribution of relative $L^2$-error for each method. The boxes represent the interquartile range, with the line inside the box showing the median error. PINO exhibits the widest box, indicating a greater spread in error values, while VINO and VINO+data display narrower boxes, suggesting lower and more consistent errors. Outliers are visible as circles outside the whiskers, highlighting points where the error was significantly different from the majority of the data. The box plot visually confirms the improvement in the performance of VINO in comparison with PINO. We note that, even by incorporating data into PINO, its accuracy is lower than VINO without data-driven training. 

\subsection{Practical Numerical Results}

To demonstrate VINO's practical applicability and efficiency in addressing complex physical scenarios, we present three representative examples in this section: porous material structure, hyperelasticity, and a plate with arbitrary holes. Each example addresses distinct challenges, including material heterogeneity, large deformation, and geometric complexity, providing an evaluation of VINO's versatility and effectiveness in real-world engineering contexts. To take advantage of optimized linear algebra operations, the JAX \cite{jax2018github} framework has been used in the following experiments.

\subsubsection{Porous Material Structure}

Porous materials are crucial in different fields like aerospace and automotive due to their ability to reduce thermal stress and maintain lightweight design resulting in unique characteristics that are not possible with homogeneous materials \cite{li2020review, audouard2024resistance}. On the other hand, designing porous structures with desirable properties requires a lot of analytical and computational efforts and numerous studies in the literature explore the mechanical analysis of porous beams \cite{babaei2022functionally, ramteke2023computational, chen2016free, kiarasi2021review, chen2023functionally, agarwal2006large}. Therefore, the goal in this example is to predict their behavior based on the porosity distribution of material and distributed load.

\begin{figure}[htbp]
    \centering
    \includegraphics[width=1.0\textwidth]{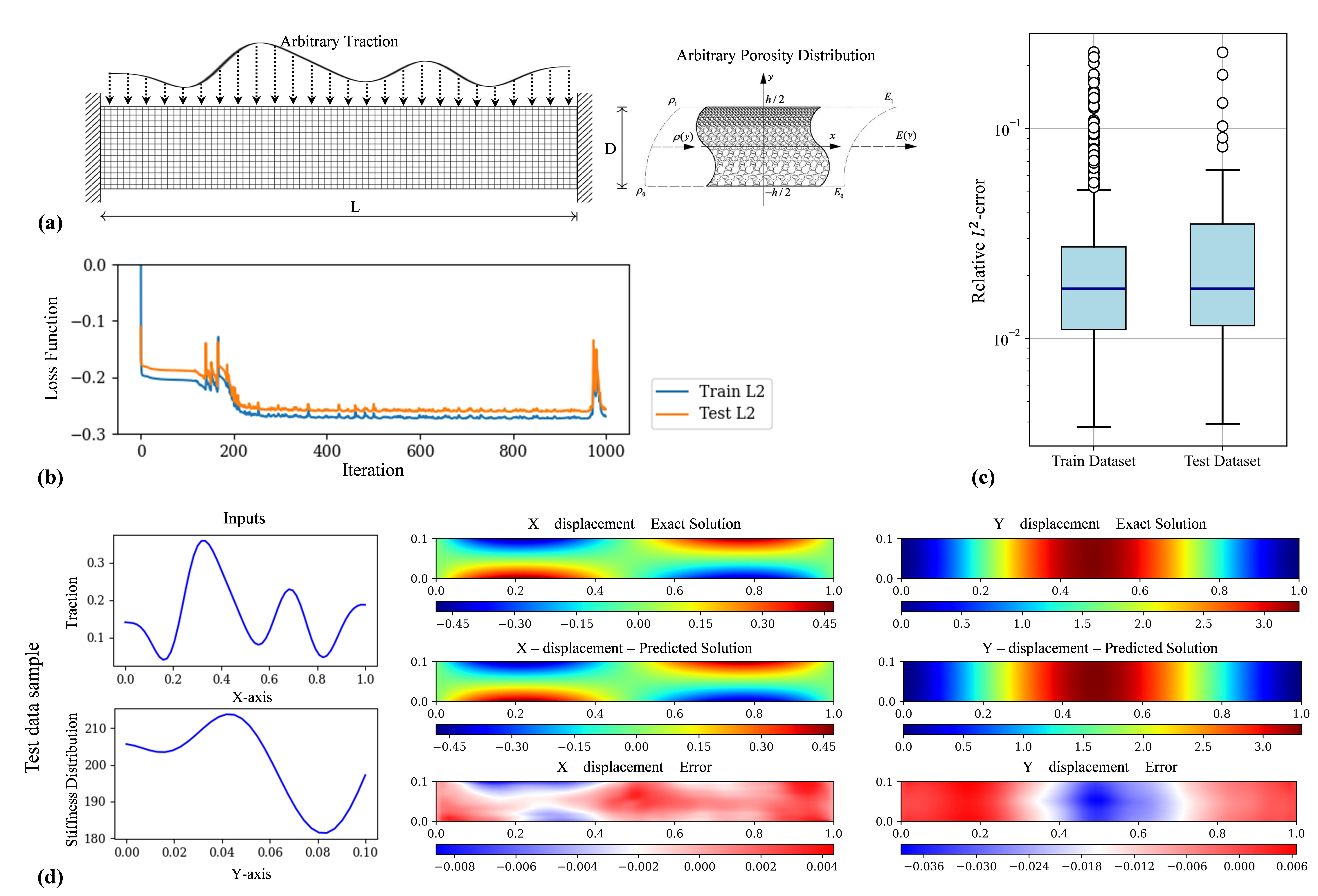}
    \caption{\textbf{Porous Material Structure:}, \textbf{(a)} Beam with Random Distributed Load and Random Porosity Distribution. \textbf{b} Loss function convergence over 1000 iterations for the train and test datasets, where both curves exhibit a rapid initial decrease followed by stabilization at lower values. \textbf{c} Box plots showing the relative \(L^2\)- error distribution for the train and test datasets, with mean values and standard deviations. The train dataset has a mean error of 0.035 with a standard deviation of 0.087, while the test dataset has a mean error of 0.032 with a standard deviation of 0.051. \textbf{(d)} Performance of VINO for predicting the displacement field of the porous material beam in a random sample from train and test dataset.}
    \label{fig: FGM}
\end{figure}

To exemplify this, VINO has been trained to learn an operator, that maps from arbitrary traction and elasticity module distribution functions into the displacement field of a beam with depth \(D\) and length \(L\). In other words, the networks are employed to solve a range of problems with different arbitrary traction functions and porosity distributions, shown in Fig. \ref{fig: FGM}a. The beam is considered to be attached to supports at \(x=0\) and \(x=L\) and sustains a distributed load on the top edge of the beam. Therefore, the governing PDE becomes as follows:

\begin{equation}
    \begin{split}
        & \frac{\partial \sigma_{xx}}{\partial x} + \frac{\partial \sigma_{xy}}{\partial y} + f_x = 0 \\
        & \frac{\partial \sigma_{yx}}{\partial x} + \frac{\partial \sigma_{yy}}{\partial y} + f_y = 0 \\        
        & \mathbf{u}(x,y) = 0 \:\:\:\quad \text{for } x=0, \, L \\
        & \sigma \cdot n = \hat{\mathbf{t}}\,(x) \:\quad  \text{for } y=D
        \label{eq: FGM Beam}
    \end{split}
\end{equation}
where \(f_x\) and \(f_y\) are body forces per unit volume and the stress components are derived based on the following constitutive relations:
\begin{equation}
\boldsymbol{\sigma} = E \mathbf{D} \boldsymbol{\varepsilon}
\end{equation}
where \(E\) is the elasticity modulus and can vary over the entire domain, matrix \(\mathbf{D}\), strain vector \(\boldsymbol{\varepsilon}\) and stress vector \(\boldsymbol{\sigma}\) are 
\begin{align}
\mathbf{D}
=
\frac{1}{1 - \nu^2}
\begin{bmatrix}
1 & \nu & 0 \\
\nu & 1 & 0 \\
0 & 0 & \frac{1 - \nu}{2}
\end{bmatrix}
\end{align}
\begin{align}
\boldsymbol{\sigma} = [\sigma_{xx},\, \sigma_{yy},\, \sigma_{xy}]^T, \quad \boldsymbol{\varepsilon} = [\varepsilon_{xx},\, \varepsilon_{yy},\, \varepsilon_{xy}]^T 
\end{align}

Here, the strain components are derived by
\begin{equation}
\varepsilon_{xx} = \frac{\partial u}{\partial x}
\end{equation}
\begin{equation}
\varepsilon_{yy} = \frac{\partial v}{\partial y}
\end{equation}
\begin{equation}
\varepsilon_{xy} = \varepsilon_{yx} = \frac{1}{2} \left( \frac{\partial u}{\partial y} + \frac{\partial v}{\partial x} \right)
\end{equation}

The desired operator from this PDE is \(\mathcal{G} \coloneqq \mathcal{L} \left\{\hat{\mathbf{t}}, \hat{E}\right\} : \mathcal{U}^t \times \mathcal{U}^E \rightarrow \mathcal{U}\) defined to map the traction and porosity distribution to the corresponding displacement field \( \left\{ \hat{\mathbf{t}}, \hat{E} \right\} \mapsto \mathbf{u}\), where \(\mathcal{U}^t\) is the space of continous real-valued functions defined on the top boundary, \(\mathcal{U}^E\) is the space of continuous real-valued functions defined on \([0,L]\times[0,D]\), and \(\mathcal{U}\) is the space of continuous functions with values in \(\mathbb{R}^2\), representing the \(x\)- and \(y\)-displacements. We solve the problem when the domain, \(\Omega\) is a rectangle with corners at \((0,0)\) and \((1,0.1)\), and \(\nu = 1/3\) is the Poisson ratio. For the data generation and validation, the isogeometric analysis (IGA) has been used \cite{anitescu2018recovery}. To create a dataset including the solution \( \left\{ \hat{\mathbf{t}}^{(i)}, \hat{E}^{(j)}, \mathbf{u}^{(i)} \right\}\), we have used a Gaussian Random Fields (GRF) \cite{williams2006gaussian} with a length scale \(l=0.1\) for traction and a length scale \(l=0.025\) for porosity distribution as 
\begin{equation}
    \begin{split}
        &\hat{\mathbf{t}}(x) \sim \sigma_{t} \times \mathcal{GP}(0, k_1(x_1, x_2)) + \mu_{t} \\
        &\hat{\mathbf{E}}(x) \sim \frac{E_{\text{max}}-E_{\text{min}}}{\mathcal{GP}_{\text{max}}-\mathcal{GP}_{\text{min}}} \times \left( \mathcal{GP}(0, k_1(x_1, x_2)) - \mathcal{GP}_{\text{min}} \right)+E_{\text{min}}
        \label{eq: FGM-Random Fields}
    \end{split}
\end{equation}
where \(\sigma_{t}\) and \(\mu_t\) are expected variance and the average of traction (here \(\sigma_{t} = 0.15\) and \(\mu_t = 0.2\)), \(E_{\text{min}}\) and \(E_{\text{max}}\) are expected values for minimum and maximum of elasticity modulus, which in the current example are selected randomly between \SI{20}{\giga\pascal} and \SI{380}{\giga\pascal}. Therefore our goal is to construct an approximation of \( \mathbf{u}^{(i)} = \mathcal{G} \left( \hat{\mathbf{t}}^{(i)}, \hat{E}^{(i)} \right) \) by the parametric map \(\mathcal{G}_{\theta} : \mathcal{U}^t \times \mathcal{U}^E \rightarrow \mathcal{U}, \, \theta \in \mathbb{R}^p\). 

The displacement field was computed using \(n_x \times n_y\) uniformly spaced points within the domain, with \(n_x = 64\) and \(n_y = 32\). Throughout all layers, the GELU (Gaussian Error Linear Unit) activation function was employed \cite{hendrycks2016gaussian}. The optimization process utilized the Adam optimizer, and the network underwent training on 1000 sets of GRF functions and was subsequently tested on 100 additional sets of GRF functions. The computational time taken for the training part is 376 s, but once the training is complete, solving the problem with any arbitrary function takes approximately 3 milliseconds, compared to 1.5 seconds for IGA. 

Figs. \ref{fig: FGM}b and \ref{fig: FGM}c display the convergence of the loss function, and the box plot, which illustrates the distribution of relative $L^2$-error for training and test dataset. Moreover, Figure \ref{fig: FGM}d shows the result of the VINO in predicting the displacement field of a sample from training and a sample from the test dataset. The VINO result is in agreement with the ground truth solution. 

\subsubsection{Hyperelasticity}

\begin{figure}[htbp]
    \centering
    \includegraphics[width=1.0\textwidth]{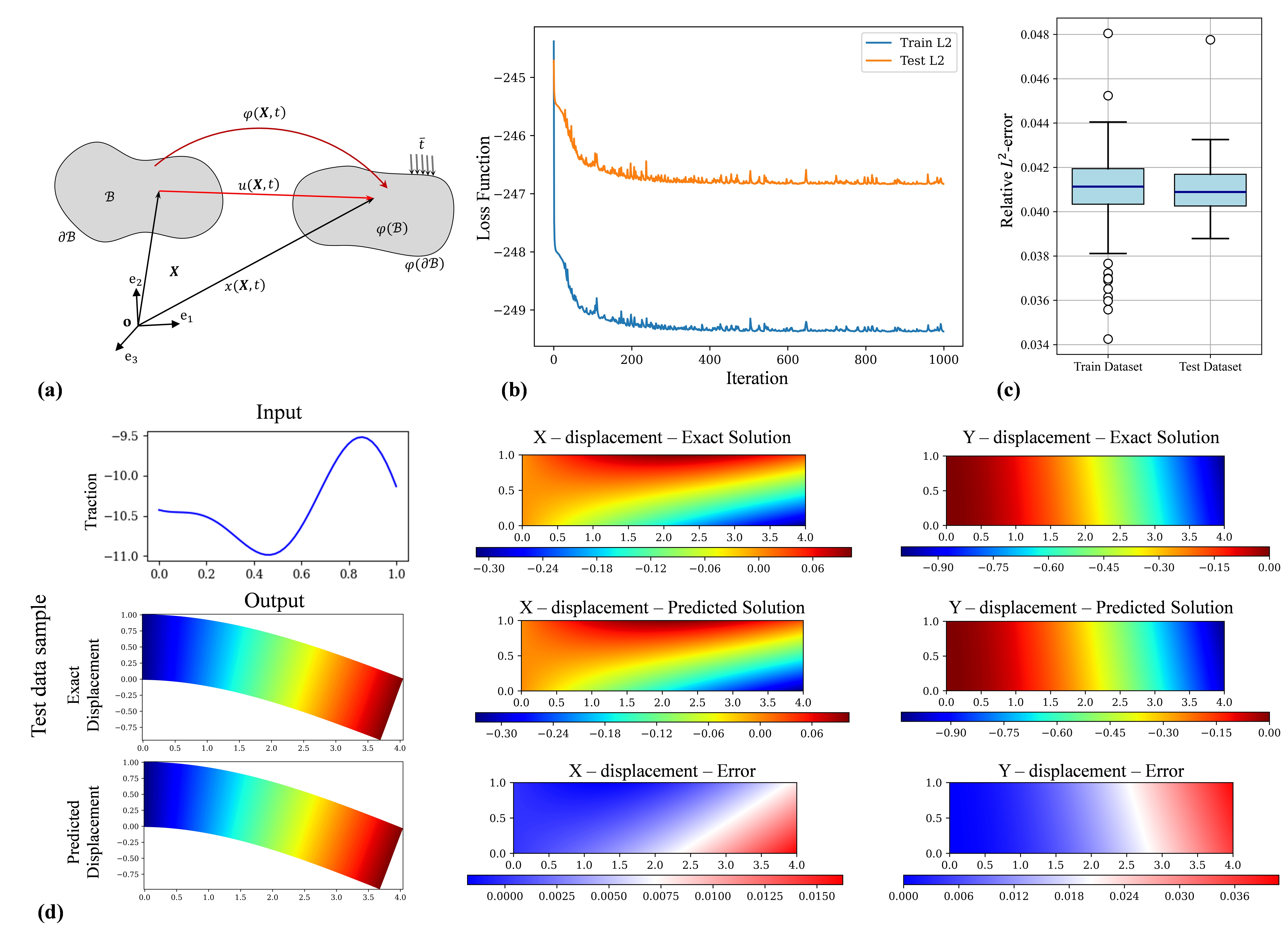}
    \caption{\textbf{Hyperelasticity.} \textbf{(a)} Motion of body \(\mathcal{B}\). \textbf{(b)} Training/testing loss convergence over 1000 iterations. The Line plot shows the model’s loss reduction and stabilization across iterations \textbf{(c)} Comparison of relative \(L^2\)-errors of VINO for train and test datasets. The box plot highlights error consistency between datasets  \textbf{(d)} Performance of VINO for predicting the displacement field of the hyperelasticity problem in a random sample from test dataset.}
    \label{fig:HyperElasticity}
\end{figure}

The next practical example is nonlinear large deformation hyperelasticity problems. Consider a body \( \mathcal{B} \) composed of a homogeneous, isotropic, nonlinear hyperelastic material. The body is represented by a set of particles or material points \( X \), bounded by the initial boundary \( \partial \mathcal{B} \). When a load \( t \) is applied over certain regions of this boundary (illustrated in Fig \ref{fig:HyperElasticity}a), the body undergoes deformation. The mapping of each material point from its initial configuration to the deformed or current configuration is given by \(\boldsymbol{\varphi}: \mathcal{B} \rightarrow \mathcal{B}_t \), where \( \boldsymbol{X} \rightarrow \boldsymbol{\varphi}(\boldsymbol{X}, t) = \boldsymbol{x} = u + \boldsymbol{X} \).

We define a boundary value problem in the initial configuration, which includes the following components:

\begin{equation}
    \begin{split}
        & \nabla \cdot \mathbf{P} + \mathbf{f}_b = \mathbf{0}  \\ 
        & \mathbf{u} = \bar{\mathbf{u}}  \:\quad \text{ on } \partial \mathcal{B}_u \\
        & \mathbf{P} \cdot \mathbf{n} = \mathbf{t}  \:\quad \text{ on } \partial \mathcal{B}_t
        \label{eq: hyperelasticity Problem}
    \end{split}
\end{equation}
where \( \nabla \cdot \mathbf{P} \) represents the divergence of the 1st Piola-Kirchhoff stress tensor \( \boldsymbol{P} \) with respect to \( \mathbf{X} \), and \( \mathbf{f}_b \) is the body force. \( \bar{\mathbf{u}} \) is the prescribed displacement on the Dirichlet boundary \( \partial \mathcal{B}_u \) and \( \mathbf{t} \) is the traction on the Neumann boundary \( \partial \mathcal{B}_t \), and \( \mathbf{n} \) is the outward unit normal vector. The boundaries are required to satisfy \( \partial \mathcal{B}_u \cup \partial \mathcal{B}_t = \partial \mathcal{B} \) and \( \partial \mathcal{B}_u \cap \partial \mathcal{B}_t = \emptyset \). The 1st Piola-Kirchhoff stress tensor \( \mathbf{P} \) is connected to the deformation gradient \( \mathbf{F} \), its power conjugate, defined by
\begin{equation}
    \mathbf{F} = \nabla \, \boldsymbol{\varphi}(\mathbf{X})
\end{equation}
through a constitutive relation, specified by
\begin{equation}
    \mathbf{P} = \frac{\partial \Psi}{\partial \mathbf{F}}
\end{equation}
where \( \Psi \) represents the strain energy density function of the material. Considering the procedure is the same for other hyperelastic models, we examined  Mooney-Rivlin hyperelastic models for \( \Psi \) 
\begin{equation}
    \Psi(I_1, I_2, J) = c(J - 1)^2 - d \log(J) + c_1 (I_1 - 3) + c_2 (I_2 - 3)
\end{equation}
where \( I_1 = \text{trace}(\mathbf{C}) \), \( I_2 = \frac{1}{2}(\text{trace}(\mathbf{C})^2 - \text{trace}(\mathbf{C} \cdot \mathbf{C})) \), and \( J = \sqrt{\det(\mathbf{C})} \). The constants \( c \), \( c_1 \), \( c_2 \), and \( d \) (where \( d = 2(c_1 + 2c_2) \)) characterize the material's properties, assuming a stress-free reference configuration and \(\mathbf{C} = \mathbf{F}^T \cdot \mathbf{F}\) is the right Cauchy-Green tensor.

In finite element analysis, solving the balance equation with the appropriate boundary conditions involves transforming the problem’s strong form to a weak form. Using the principle of virtual displacements, linear momentum can be derived in the initial configuration. However, there exists a strain energy, which is the elastic energy stored in the body \( \Psi \) in hyperelastic material. The classical principle of the minimum potential energy can be reformulated based on this strain energy. In the case of elastostatics, the potential is given by:
\begin{equation}
    \Pi(\boldsymbol{\varphi}) = \int_\mathcal{B} \Psi \, dV - \int_\mathcal{B} \mathbf{f}_b \cdot \boldsymbol{\varphi} \, dV - \int_{\partial \mathcal{B}_t} \mathbf{t} \cdot \boldsymbol{\varphi} \, dA
\end{equation}
To find a stationary point, we minimize the potential energy:
\begin{equation}
    \min_{\boldsymbol{\varphi} \in H} \Pi(\boldsymbol{\varphi})
\end{equation}
after discretization of domain, where \( H \) is the space of admissible (trial) functions. For analytical calculation of the integral inside each element, the fourth-order Taylor series has been used instead of \(\log\) function, taking advantage of its simplicity. 

To exemplify this problem, we considered plane strain condition in a beam problem with length \(L = \SI{4}{m}\) and height \(H = \SI{1}{m}\) and arbitrary traction load at the right end and is clamped at the left side.  The goal is mapping from arbitrary traction into displacement fields. In other words, the networks are employed to solve a range of problems with different arbitrary traction functions. The ground truth solution is obtained by FEniCS \cite{BarattaEtal2023} on a fine mesh.  The displacement field was computed using \(n_x \times n_y\) uniformly spaced points within the domain, with \(n_x = 200\) and \(n_y = 50\). The optimization process utilized the Adam optimizer, and the network underwent training on 500 sets of GRF functions and was subsequently tested on 50 additional sets of GRF functions. 

Figs. \ref{fig:HyperElasticity}b and \ref{fig:HyperElasticity}c  display the convergence of the loss function, beside the box plot, which compares the distribution of the relative $L^2$-error for training and test dataset. The median error is represented by the blue line within each box, while the box bounds represent the interquartile range. This comparison helps assess the consistency of model performance on both datasets. Moreover, the right plot shows the training and testing loss values across 1000 iterations. The blue line represents the training loss, and the orange line represents the testing loss. The initial steep decline indicates that the model quickly reduces its loss, followed by convergence to a steady state as training progresses. The stability in both curves at later iterations indicates that the model has reached a stable level of performance with minimal fluctuations. In addition, Figure \ref{fig:HyperElasticity}d shows the result of the VINO in predicting the displacement field of a sample from the test dataset. The VINO result has a very good agreement with the ground truth solution. The computational time taken for the training part is 1282 s, but once the training is complete, solving the problem with any arbitrary function takes approximately 0.07 s, compared to 40 seconds on average for FEniCS. 

\subsubsection{Plate with arbitrary voids}

\begin{figure}[htbp]
    \centering
    \includegraphics[width=1\textwidth]{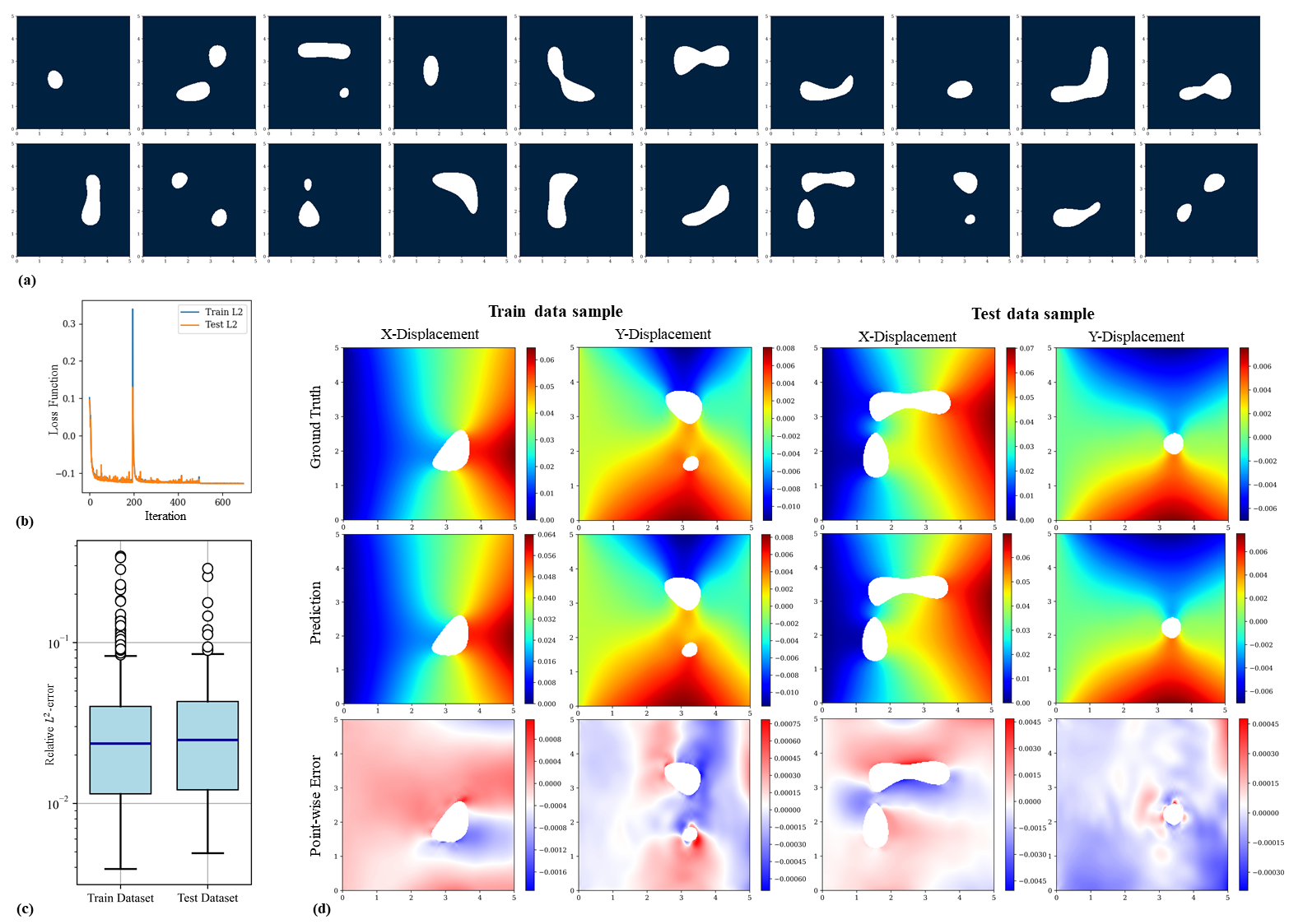}
    \caption{\textbf{Plate with arbitrary voids.} \textbf{(a)} Random samples of allowable domain with arbitrary void shapes. \textbf{(b)} Loss function convergence over 1000 iterations for the train and test datasets, where both curves exhibit a rapid initial decrease followed by stabilization at lower values. \textbf{(c)} Box plots showing the relative \(L^2\)- error distribution for the train and test datasets, with mean values highlighted in blue and standard deviations in red. The train dataset has a mean error of 0.035 with a standard deviation of 0.087, while the test dataset has a mean error of 0.032 with a standard deviation of 0.051. \textbf{(d)} Performance of VINO for predicting the displacement field of the porous material beam in a random sample from train and test dataset}
    \label{fig:Void}
\end{figure}

As a third practical example, we consider parametric partial differential equations like Eq. \ref{eq: FGM Beam} but on various domains. Assume the problem domain \(\Omega_p\) is parameterized by design parameters \(p \in \mathcal{P}\), which follows some distribution \(p \sim \mu\). Therefore, the problem setting becomes:
\begin{equation}
    \begin{split}
        \frac{\partial \sigma_{xx}}{\partial x} + \frac{\partial \sigma_{xy}}{\partial y} + f_x = 0  & \:\quad\quad \text{for } x, y \in \Omega_p\\
        \frac{\partial \sigma_{yx}}{\partial x} + \frac{\partial \sigma_{yy}}{\partial y} + f_y = 0  & \:\quad\quad \text{for } x, y \in \Omega_p \\
        \mathbf{u}(x,y) = \mathbf{\hat{u}} & \:\quad\quad \text{for } x=0\\
        \boldsymbol{\sigma} \cdot \boldsymbol{n} = \mathbf{\hat{t}}\,(y) &  \:\quad\quad  \text{for } x=L
        \label{eq: arbitrary hole}
    \end{split}
\end{equation}
where \(\mathbf{\hat{u}} \in \mathcal{U}\) is Dirichlet boundary condition, the \(\mathbf{\hat{u}} \in \mathcal{T}\) is the Neumann boundary condition and the goal is approximation mapping \(\mathcal{G}\) from problem domain to the PDE solution, \(\Omega_p \mapsto u\). In this example, we assumed Dirichlet and Neumann boundary conditions are fixed and the domain space is given as meshes (uniform grid) \(p^{(i)} = \{ x^{(I)}, y^{(I)} \} \subset \Omega_p\). We consider \(\Omega_p = [0,5]^2 -\mathcal{V} \), where \(\mathcal{V}\) is composed of arbitrary shaped holes inside the region \([1,4]^2\). For more clarification, some samples are depicted in Figure \ref{fig:Void}a. The plate is clamped on the left edge and constant tension traction is applied on the right edge. 

We generated 1000 training data and 200 test data and used an IGA solver for the ground truth solution. Note that the ground truth solutions are not used during the training process and are just for comparing the final result of the network. The displacement field was computed using \(n_x \times n_y\) uniformly spaced points within the domain, with \(n_x = n_y = 200\) and the optimization process utilized the Adam optimizer.  Figures \ref{fig:Void}b and \ref{fig:Void}c illustrate both the distribution of the relative  $ L^2$ error and the convergence of the loss function. The box plot the relative $ L^2$ error for the training and test datasets, with the blue line inside each box representing the median error. The boxes capture the interquartile range, allowing for an evaluation of the model's error consistency across both datasets. The line plot depicts the training and testing loss values over 1000 iterations. The blue line corresponds to the training loss, and the orange line to the testing loss. The sharp initial drop indicates a rapid reduction in loss, followed by stabilization, suggesting that the model achieves a steady state with minimal fluctuation as training progresses. In addition, Figure \ref{fig:Void} shows the result of the VINO in predicting the displacement field of some samples from the training and test dataset. The VINO result has an acceptable agreement with the ground truth solution.

\section{Discussion} \label{sec:Discussion}
In this work, we introduced the Variational Physics-Informed Neural Operator (VINO) that overcomes the limitations of data-driven methods that rely on large datasets of parid input and output, which are often costly or unavailable. VINO's core innovation lies in its unique incorporation of the variational form of physical laws directly into the learning process by discretizing the domain and using shape functions to analytically handle derivatives and integrals. This approach allows VINO to learn operators solely from physical principles, a capability that positions it as a unique tool among neural operators. Our results demonstrate that VINO can outperform existing approaches like PINO, particularly as mesh resolution increases. Notably, while other methods exhibit decreasing accuracy or inconsistent performance as mesh size grows, VINO demonstrates convergent behavior in fine meshes, highlighting its robustness and scalability. In final, VINO suggests a better way to incorporate the physical laws into neural operators.

Currently, VINO is configured to solve PDEs with an available variational form. Future research should extend this approach to the weak form of equations, making it applicable to all types of PDEs. This study implemented VINO for rectangular elements, but the same methodology could be adapted for other element types, such as triangular or tetrahedron elements. Additionally, future work will refine VINO’s architecture to handle more complex, multi-scale, and non-linear problems, explore adaptive meshing strategies, and expand applications to domains where traditional solvers encounter computational limits. VINO thus marks a promising advancement in neural operators for PDE learning, blending physics with machine learning to enhance operator learning frameworks. The development of a software library of pre-trained models also holds significant potential.

\section{Methods} \label{sec: Method}
In this section, we introduce the PDE systems under consideration. We then provide an overview of the operator-learning and physics-informed frameworks. Finally, we present our approach in detail in the next subsections. 

\subsection{Problem setting} \label{sec:Preliminaries}
If \((\mathcal{V}, \mathcal{S}, \mathcal{B})\) is a triplet of Banach spaces, and \( \mathcal{N}: \mathcal{V} \times \mathcal{S} \to \mathcal{B} \) is a linear or nonlinear differential operator, we consider parametric PDEs taking the form 

\begin{equation}
    \begin{split}        
        \mathcal{N}(p, s) &= 0, \quad \text{in } \Omega \subset \mathbb{R}^d \\
        s &= \bar{s}, \quad \text{in } \partial \Omega_D \\
        \mathcal{L} (s) &= \bar{t}, \quad \text{in } \partial \Omega_N
        \label{eq: PDE}
    \end{split}
\end{equation}

where \(p \in \mathcal{P} \subset \mathcal{V}\) is the PDE parameters (i.e. coefficient or input functions), \(s \in \mathcal{S} \) is the corresponding unknown solutions, \(\Omega\) is a bounded domain and \(\mathcal{L}\) is a differential operator. Typically, the functions \(\bar{s}\) and \(\bar{t}\) serve as fixed boundary conditions, but they may also be included as a parameter.
Then, we can consider the solution operator \( \mathcal{G}: \mathcal{P} \to \mathcal{S}\) which maps the parameters to the solutions \( p \mapsto s \). 

\subsubsection{Solving equation in strong form using PINNs}
This section discusses the problem of approximating a solution \(s^\dagger\) for a given instance \(p^\dagger\) using a solution operator \(\mathcal{G}^\dagger\) defined by Eq. \ref{eq: PDE}. The task involves finding an approximation to the true solution by implementing ML-enhanced conventional solvers \cite{kochkov2021machine, pathak2021mlP, greenfeld2019learning} as well as neural network-based solvers such as the PINNs and Deep Galerkin Method \cite{raissi2019physics, mishra2024artificial, sirignano2018dgm}. These PINN-type approaches use a neural network \( s_\theta \) with parameters \( \theta \) to approximate the solution function \( s^\dagger \). The parameters \( \theta \) are derived by minimizing the loss function based on the PDE strong form with exact derivatives computed utilizing the AD method. In other words, the physics-informed loss is defined by minimizing the left-hand side of Eq. \ref{eq: PDE} in the squared norm of \( \mathcal{B} \), giving the loss function:
\begin{equation}
    \begin{split}        
        \mathcal{L}_{\text{pde}}(p^\dagger, s_\theta) & = \alpha \mathcal{L}_{\text{system}} + \beta \mathcal{L}_{\text{BC}}
        \\
        & = \alpha \left\|\mathcal{N}(p^\dagger, s_\theta)\right\|_{L^2(\Omega)}^2 + \beta \left\| s_\theta |_{\partial \Omega} - \bar{s} \right\|_{L^2(\partial \Omega)}^2,
        \label{eq: PINN_loss}
    \end{split}
\end{equation}
where \(\alpha\), \(\beta > 0\) are hyperparameters. In this paper, the notation \(\dagger\) indicates that the corresponding quantity is fixed and specified.

PINNs offer the advantage of incorporating physics directly into the learning process and the universal approximability of neural networks but, in practice, PINNs often struggle to solve complex PDEs, particularly when the solution contains multi-scale features of high-frequency terms \cite{wang2022and}, or nonlinear relationships of hyperbolic PDEs, such as the problem of immiscible two-phase fluid transport in porous media \cite{fuks2020limitations}, particularly in regions with shockwaves. whereas diffusion term in parabolic PDEs ensures improved data estimation \cite{huang2023limitations}. Moreover, in dealing with large-scale problems, PINNs require efficient implementations of multi-GPU architectures, which is not always available in practice, and advanced parallelization strategies—such as data parallel or hybrid data and model parallel paradigms \cite{cai2021physics}. Additionally, PINNs solve PDEs for only a single set of parameters at a time and when the parameters change, the model requires retraining, further increasing its inefficiency \cite{sadegh2024deepnetbeam}. In this study, we aim to improve efficiency by combining neural operators and PINNs in a variational format.  

\subsubsection{Solving the variational form using the physics-informed neural networks}

For problems in which a variational principle exists, a functional \(\Pi\) can be derived in the form of 
\begin{equation}
    \begin{split}        
        \Pi (p, s) = \int_\Omega \mathcal{F}\left(p, s, \frac{ds}{dx}, \ldots\right) dx + \int_{\partial\Omega} \mathcal{E}\left(p, s, \frac{ds}{dx}, \ldots\right) dx
        \label{eq: variational_principle}
    \end{split}
\end{equation}
where \(\mathcal{F}\) and  \(\mathcal{E}\) are specified differential operators and the solution to the continuous problem is a function \(s\) which makes \(\Pi\) stationary. For instance, this concept is equivalent to the principle of least work in mechanics, which posits that the solution to a static mechanics problem results in the total energy—defined as the strain energy minus the external work performed by external forces at the boundary and body forces—reaching a minimum. Additionally, this principle mirrors Fermat's principle in optics (also known as the principle of least time), which asserts that the path taken by a ray of light between two points is the one that can be traversed in the shortest time. For example, the loss function in mechanical problems is formulated as \cite{samaniego2020energy} 
\begin{equation}
    \begin{split}        
        \mathcal{L}_{\text{pde}}(p^\dagger, s) = \frac{1}{2} \int_\Omega \varepsilon(u_\theta) : \mathbb{C} : \varepsilon(u_\theta) \, d\Omega - \left( \int_{\Omega} f \cdot u_\theta \, d\Omega + \int_{\partial \Omega_N} \bar{t} \cdot u_\theta \, d\Gamma \right), 
        \label{eq: VPINN_loss}
    \end{split}
\end{equation}
where \(\varepsilon\) denotes strain tensor, which in turn depends on the displacement \(u_\theta\) predicted by a neural network with parameters \(\theta\) as solution \(s\), and constitutive elastic matrix \(\mathbb{C}\), prescribed boundary forces \(\bar{t}\) and body force \(f\) serve as the fixed PDE parameters \(p^\dagger\).

Similar to the PINN loss function in strong form, Eq. \ref{eq: PINN_loss}, the VPINN loss function, Eq. \ref{eq: VPINN_loss} is also defined for a set of specific PDE parameters. However, there are two key differences: 

\begin{itemize}
    \item The PINN loss function, in essence, represents a set of equations that the networks attempt to satisfy at the collocation points, with the exact solution satisfying all these equations simultaneously. In contrast, the VPINN loss function is not tied to direct equations; instead, it is a functional that the exact solution minimizes. The closer the network's predicted solution is to minimizing this functional, the closer it is to the exact solution.
    \item While the PINN loss function is designed to reduce the squared error (i.e., the norm) at the collocation points, with a lower limit of zero, the VPINN loss function has no defined lower bound and can take on negative values. In other words, in physical problems, the PINN loss function represents an error and can be dimensionless. Conversely, the VPINN loss function is not an error; it has physical dimensions and meaning. For example, in mechanical problems, it represents potential energy.
\end{itemize}

VPINN necessitates fewer hyperparameters compared to the strong form because it incorporates the Neumann boundary condition in the physical part of the loss function and does not need to define tuning parameters in the loss function. In addition, VPINN achieves greater computational efficiency and accuracy due to its lower derivative order \cite{wang2022cenn}. However, VPINN is dependent on the chosen integration scheme for loss calculation, which must be accurate for effective numerical integration \cite{li2021physics, sheng2021pfnn}, and it lacks generality since not all PDEs possess a corresponding variational form \cite{zienkiewicz2005finite}. Current research primarily concentrates on the strong form of PINN, with relatively little focus on the variational form. The strong form involves numerous hyperparameters, particularly a penalty factor or the multiplication of coordinates to meet specific geometric boundary conditions \cite{jagtap2020conservative}.

\subsubsection{Learning the Solution Operator Via Neural Operators}

An alternative approach is to learn the solution operator \( \mathcal{G} \). Given a PDE described by equation \ref{eq: PDE} and its corresponding solution operator \( \mathcal{G}\), a neural operator \( \mathcal{G}_\theta \) with parameters \( \theta \) can be used as a surrogate model to approximate \( \mathcal{G} \). Typically, we assume that a data set \( \{p_j, s_j\}_{j=1}^N \) is available, where \( \mathcal{G}(p_j) = s_j \) and \( p_j \sim \mu \) are independent and identically distributed (i.i.d.) samples from a distribution \( \mu \) supported on \( \mathcal{P} \). In this scenario, optimizing the solution operator involves minimizing the empirical data loss for a given data pair:

\begin{equation}
    \begin{split}        
        \mathcal{L}_{\text{data}}(s, \mathcal{G}_\theta(p)) = \left\|s - \mathcal{G}_\theta(p)\right\|_\mathcal{S}^2 =  \int_{\Omega} |s(x) - \mathcal{G}_\theta(p)(x)|^2 dx  
        \label{eq:Neural_Operator_loss}
    \end{split}
\end{equation}

The neural operator has been developed as an extension of conventional deep neural networks, to facilitate the handling of operators. The approximation of nonlinear operators is achieved by composing a linear integral operator, denoted by \(\mathcal{K}\), with a pointwise nonlinear activation function, denoted by \(\sigma\). The neural operator, designated as \(\mathcal{G}_\theta\), is defined as follows:
\begin{equation}
    \begin{split}        
        \mathcal{G}_\theta :=\mathcal{Q} \circ (\mathcal{W}_L + \mathcal{K}_L) \circ \cdots \circ \sigma(\mathcal{W}_1 + \mathcal{K}_1) \circ \mathcal{P},
        \label{eq:Neural_Operator}
    \end{split}
\end{equation}
wherein  \(\mathcal{P}\) and  \(\mathcal{Q}\) are pointwise operators parameterized by neural networks. In particular, the operator \(\mathcal{P}\) maps from the space of dimension \(d_p\) to the space of dimension \(d_1\), while the operator \(\mathcal{Q}\) maps from the space of dimension \(d_L\) to the space of dimension \(d_u\). The operator \(\mathcal{P}\) elevates the lower-dimensional input function to a higher-dimensional space, whereas \(\mathcal{Q}\) reduces the high-dimensional output function to a lower-dimensional space.

The model stacks \(L\) layers, each of the form \(\sigma(\mathcal{W}_l + \mathcal{K}_l)\), where \(\mathcal{W}\) are pointwise linear operators represented as matrices \(\mathcal{W}_l \in \mathbb{R}^{d_{l+1} \times d_l}\) and the kernel operators, denoted by \(\mathcal{K}_l: \{\Omega \to \mathbb{R}^{d_l}\} \to \{\Omega \to \mathbb{R}^{d_{l+1}}\}\), are integral operators. The activation functions \(\sigma\) are fixed and the parameters, \(\theta\) include all the parameters of  \(\mathcal{P}\),  \(\mathcal{Q}\),  \(\mathcal{W}_l\), and \(\mathcal{K}_l\). For more details, we refer to \cite{li2020fourier}.

As with supervised learning in computer vision and language, operator learning relies heavily on data, necessitating the availability of training and testing points that adhere to the same distribution. However, neural operator models like FNO have difficulty generalizing across different coefficients or geometries. Furthermore, gathering representative data is particularly challenging for complex PDEs with slow or non-existent solvers. Additionally, FNO models do not incorporate equation knowledge, which prevents them from fully closing the generalization gap and limits their broader application \cite{Goswami2023}. The methods of the physics-informed neural operator (PINO) \cite{li2024physics} and physics-informed DeepOnets \cite{wang2021learning} are introduced to overcome these challenges, but the reliance on data and associated issues still persist. We introduce the VINO method to address these problems in the following methods and demonstrate that dependence on data can be overcome.

\subsection{Learning the Solution Operator Via Variational Physics-Informed Neural Operator (VINO)}

The limitation of data-dependence approaches, including data-driven and physics-informed neural operators, is that the solution requires a set of paired datasets, which can be challenging to gather or generate, particularly in engineering applications where data acquisition is often prohibitively expensive. With motivation from Deep Ritz \cite{yu2018deep} and DEM \cite{samaniego2020energy} we used the functional minimization in the variational format of PDE to train the network, obviating the necessity for a paired dataset. Although using the variational form as a loss function is mentioned in the PINO paper for solving the Darcy equation, the authors abandoned it, citing the superior performance of the strong form. However, the variational form allows for domain discretization, which we inspired by the Finite Element Method, resulting in a notable enhancement. The most significant improvement involves removing the necessity of a dataset and enhancing the method's robustness to resolution changes.

Consider again the functional presented in equation \ref{eq: variational_principle}. We use it directly as a loss function, so the neural network \(\mathcal{G}_{\theta}\) tries to minimize it across all training function inputs:
\begin{equation}
    \begin{split}        
        \mathcal{L}_{\text{pde}} (\theta) =\frac{1}{n} \sum_{i=1}^{n} \ \int_\Omega \mathcal{F}\left(p^{(i)}, \mathcal{G}_{\theta}(p^{(i)}), \frac{d\mathcal{G}_{\theta}(p^{(i)})}{dx}, \ldots\right) d\Omega 
        +
        \int_{\partial\Omega}
        \mathcal{E}\left(p^{(i)}, \mathcal{G}_{\theta}(p^{(i)}), \frac{d\mathcal{G}_{\theta}(p^{(i)})}{dx}, \ldots\right) d\Omega
        \label{eq: Variational_loss_function}
    \end{split}
\end{equation}

The main bottleneck in deriving this loss function is the computation of the numerical integration. Various methods, such as trapezoidal, Simpson's, and Monte Carlo, have been used in the literature, but none have proven sufficient for neural operators. Furthermore, a significant technical challenge lies in efficiently computing derivatives for neural operators \cite{li2024physics}. Numerical differentiation methods, such as finite differences, are fast and memory-efficient but encounter challenges. On one hand, they require a fine-resolution grid, while on the other, increasing the resolution increases the computational cost and may lead to numerical instabilities. To deal with these challenges, inspired by the finite element method, we can divide the domain into smaller elements and evaluate analytically the integrals individually over each element by noting 
\begin{equation}
    \begin{split}        
        \mathcal{L}_{\text{pde}} (\theta) =\frac{1}{n} \sum_{i=1}^{n}  \sum_{e}
        \int_{\Omega_e} \mathcal{F}\left(p^{(i)}, \mathcal{G}_{\theta}(p^{(i)}), \frac{d\mathcal{G}_{\theta}(p^{(i)})}{dx}, \ldots\right) d\Omega_e
        +
        \int_{\partial\Omega_e}
        \mathcal{E}\left(p^{(i)}, \mathcal{G}_{\theta}(p^{(i)}), \frac{d\mathcal{G}_{\theta}(p^{(i)})}{dx}, \ldots\right) d\Omega_e 
        \label{eq: VINO_loss_function}
    \end{split}
\end{equation}

In the following, for clarification, we derive the equations for a rectangle element with four nodes, but the process is the same for other types of elements. The rectangular element considered has side length of \(2a\) and \(2b\) in the \(x\)- direction and \(y\)- direction, respectively. It is convenient to use a local Cartesian system \(x'\), \(y'\) defined by 
\begin{equation}
    \begin{split}        
        x' = x - x_0 \quad\quad \text{and} \quad\quad y' = y - y_0
        \label{eq: local system}
    \end{split}
\end{equation}
where
\begin{equation}
    \begin{split}        
        x_0 = \frac{1}{4} \sum_{a=1}^4 x_a \quad\quad \text{and} \quad\quad y_0 = \frac{1}{4} \sum_{a=1}^4 y_a
        \label{eq: element center}
    \end{split}
\end{equation}
in which \(x_0\), \(y_0\) are located at the center of the rectangle and \(x_a\), \(y_a\) are coordinates of the nodes. To construct approximating functions \(\hat{\phi}\) for each of the solution components, denoted by \(\phi\), these functions must have linear behavior along each edge of the element to ensure interelement \(C_0\) continuity. Thus, a suitable choice is given by
\begin{equation}
    \begin{split}        
        \hat{\phi}^e = \alpha_1 + x' \alpha_2 + y' \alpha_3 + x'y' \alpha_4 
        \label{eq: element approximation function}
    \end{split}
\end{equation}

The coefficients \(\alpha_a\) may be obtained by expressing Eq. \ref{eq: element approximation function} at each node giving
\begin{equation}
\begin{bmatrix}
\tilde{\phi}_1^e \\
\tilde{\phi}_2^e \\
\tilde{\phi}_3^e \\
\tilde{\phi}_4^e
\end{bmatrix}
=
\begin{bmatrix}
1 & -a & -b &  ab \\
1 &  a & -b & -ab \\
1 &  a &  b &  ab \\
1 & -a &  b & -ab
\end{bmatrix}
\begin{Bmatrix}
\alpha_1 \\
\alpha_2 \\
\alpha_3 \\
\alpha_4
\end{Bmatrix}
\label{eq: coefficient system}
\end{equation}

We can solve for \(\alpha_a\) in terms of the nodal values to obtain finally 
\begin{equation}
    \begin{split}        
        \hat{\phi}_e & = \frac{1}{4} \left( 1 - \frac{x'}{a} \right) \left( 1 - \frac{y'}{b} \right) \tilde{\phi}_1^e
                       + \frac{1}{4} \left( 1 + \frac{x'}{a} \right) \left( 1 - \frac{y'}{b} \right) \tilde{\phi}_2^e \\
                     & + \frac{1}{4} \left( 1 + \frac{x'}{a} \right) \left( 1 + \frac{y'}{b} \right) \tilde{\phi}_3^e
                       + \frac{1}{4} \left( 1 - \frac{x'}{a} \right) \left( 1 + \frac{y'}{b} \right) \tilde{\phi}_4^e,
        \label{eq: approximative distribution}
    \end{split}
\end{equation}

From Eq. \ref{eq: approximative distribution} we can define the four shape function as follows:
\begin{equation}
    \begin{split}        
        N_1 = \frac{1}{4} \left( 1 - \frac{x'}{a} \right) \left( 1 - \frac{y'}{b} \right), \quad
        N_2 = \frac{1}{4} \left( 1 + \frac{x'}{a} \right) \left( 1 - \frac{y'}{b} \right) \\
        N_3 = \frac{1}{4} \left( 1 + \frac{x'}{a} \right) \left( 1 + \frac{y'}{b} \right), \quad
        N_4 = \frac{1}{4} \left( 1 - \frac{x'}{a} \right) \left( 1 + \frac{y'}{b} \right)
        \label{eq: shape function}
    \end{split}
\end{equation}

As the shape functions in Eq. \ref{eq: shape function} is similar for all solution components, we can write
\begin{equation}
    \begin{split}        
        \mathcal{G}_\theta^e(p) = s^e \approx \hat{s}^e = \sum_{a=1}^4 N_a \tilde{s}_a = \mathbf{N}\mathbf{\tilde{s}}
        \label{eq: approximative solution}
    \end{split}
\end{equation}

Therefore, we can calculate the derivative and integral analytically, thus solving the difficulties and challenges mentioned above. For more clarification, we derive the equation for Poisson equation in the form of 
\begin{equation}
    \begin{split}        
        -\nabla^2 \phi = f \quad \text{in} \ \Omega
        \label{eq: Poisson}
    \end{split}
\end{equation}
where \( \nabla^2 \) is the Laplace operator, \( \phi \) is the unknown function, and \( f \) is a known source term.

The equivalent minimization problem seeks to minimize the following functional:
\begin{equation}
    \begin{split}        
        \Pi(\phi) = \frac{1}{2} \int_{\Omega} |\nabla \phi|^2 \, d\Omega - \int_{\Omega} f \phi \, d\Omega
        \label{eq: Poisson Functional}
    \end{split}
\end{equation}

So the loss function, after the discretization of the domain, is defined as 
\begin{equation}
    \begin{split}        
        \mathcal{L}_{\text{pde}} (\theta) & =\frac{1}{n} \sum_{i=1}^{n}  \sum_{e}
        \int_{\Omega_e} \left(\frac{1}{2} |\nabla \phi|^2 - f \phi \right) d\Omega_e
        \label{eq: Poisson_loss_function}
    \end{split}
\end{equation}

Since \( \phi \approx N \tilde{\phi} \), where:
\begin{equation}
    N = [N_1, N_2, N_3, N_4], \quad \tilde{\phi} = [\phi_1, \phi_2, \phi_3, \phi_4]^T.
\end{equation}

Thus, the gradient of \( \phi \) is given by:
\begin{equation}
\nabla \phi = B \tilde{\phi}, 
\end{equation}
where \( B \) is the gradient matrix of the shape functions:
\begin{equation}
B = 
\begin{bmatrix}
\frac{\partial N_1}{\partial x} & \frac{\partial N_2}{\partial x} & \frac{\partial N_3}{\partial x} & \frac{\partial N_4}{\partial x} \\
\frac{\partial N_1}{\partial y} & \frac{\partial N_2}{\partial y} & \frac{\partial N_3}{\partial y} & \frac{\partial N_4}{\partial y}
\end{bmatrix}.
\end{equation}

Therefore, the magnitude squared of the gradient is:
\begin{equation}
|\nabla \phi|^2 = (\nabla \phi) \cdot (\nabla \phi) = (B \tilde{\phi})^T (B \tilde{\phi}) = 
\tilde{\phi}^T B^T B \tilde{\phi}.
\end{equation}

As a result:
\begin{equation}
\frac{1}{2} \int_{\Omega_e} |\nabla \phi|^2 \, d\Omega_e = \frac{1}{2} \tilde{\phi}^T \left( \int_{\Omega_e} B^T B \, d\Omega_e \right) \tilde{\phi} = \frac{1}{2} \tilde{\phi}^T K^e \tilde{\phi},
\end{equation}
where \( K^e = \int_{\Omega_e} B^T B \, d\Omega_e \) is a scalar quantity and the next term becomes:
\begin{equation}
\int_{\Omega_e} f \phi \, d\Omega_e = \tilde{\phi}^T \left( \int_{\Omega_e} N^T N \, d\Omega_e \right) \tilde{f} = \tilde{\phi}^T M^e \tilde{f},
\end{equation}
where \( M^e = \int_{\Omega_e} N^T N \, d\Omega_e \) is also a scalar quantity.

Thus, the final expression, after calculating the derivatives and integral, for functional \( \Pi(\phi) \) is:
\begin{equation}
\Pi(\phi) = \frac{1}{n} \sum_{i=1}^{n}  \sum_{e} \left[ \frac{1}{2} \tilde{\phi}^T K \tilde{\phi} - \tilde{\phi}^T M \tilde{f} \right].
\label{eq: final Poisson functional}
\end{equation}

Therefore, to approximate the solution operator of the Poisson equation, VINO only needs to minimize Eq. \ref{eq: final Poisson functional} without worrying about calculating the derivatives and integral since they have already been calculated analytically inside each element. Note that the result of the Eq. \ref{eq: final Poisson functional} in the whole domain is simply the summation of it over the elements and there is no concept like assembling as there is in finite element methods. The same process can be considered for other types of PDEs as well as other types of elements, which we do not mention for brevity.

To investigate the efficiency of the discretization in VINO, we compare the performance of a neural operator with the same hyperparameters and a variational loss function format. We employ common integration methods used in the literature, including the Midpoint Method, Trapezoidal Rule, and Simpson's Rule, where derivatives are computed using the Finite Differences method. These results are then compared with VINO. To emphasize the improvement offered by VINO, we use a coarse \(16 \times 16\) grid. 

\begin{figure}[tb]
    \centering
    \includegraphics[width=1.0\textwidth]{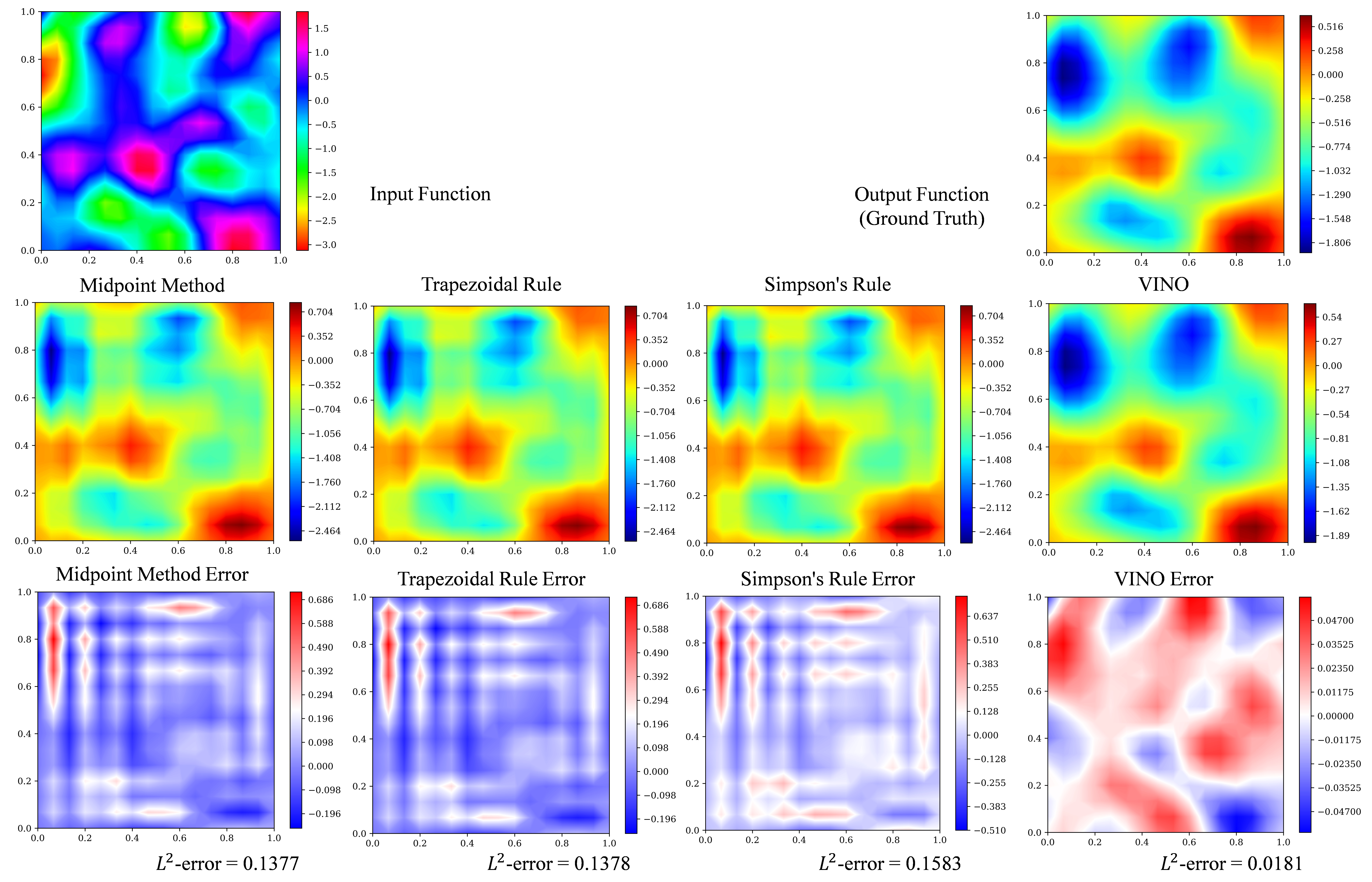}
    \vspace{-20pt}
    \caption{Comparison of the performance of different common integral methods in the literature, and VINO for Poisson Equation with a 16\(\times\)16 grid. The relative \(L^2\) errors are indicated for each model.}
    \label{fig: Figure-Integration-Methods}
\end{figure}

In Figure \ref{fig: Figure-Integration-Methods}, the top row shows the input and output functions and the middle row visualizes the predictions for all the methods. The bottom row quantifies the errors spatially, illustrating how the error is distributed across the domain. These plots show the point-wise difference between each method’s output and the ground truth solution. The corresponding \(L^2\)-error for each method is noted at the bottom of these plots \(L^2\)-error for the Midpoint Method, Trapezoidal Rule, Simpson's Rule, and VINO are 0.1377, 0.1378, 0.1583 and 0.0181 respectively.  These results demonstrate that traditional numerical integration methods and finite difference approaches for derivatives are not sufficient when using a variational loss function. In contrast, the significantly lower \(L^2\)-error for VINO highlights its superior accuracy, even on a coarse mesh.

\printbibliography

\section*{Declaration of Competing Interest}
The authors declare that they have no known competing financial interests or personal relationships that could have appeared to influence the work reported in this paper.

\section*{Acknowledgments}
The authors would like to acknowledge the support provided by the German Academic Exchange Service (DAAD) through a scholarship awarded to Mohammad Sadegh Eshaghi during this research.

\end{document}